\documentclass[12pt]{article}
\usepackage{amsfonts,amssymb,amsmath,amscd,amsthm} 

\newtheorem{thm}{Theorem}
\newtheorem{lem}[thm]{Lemma}
\newtheorem{prop}[thm]{Proposition}

\newtheorem{cor}[thm]{Corollary}
\newtheorem{prob}[thm]{Problem}
\newtheorem{conj}[thm]{Conjecture}
\theoremstyle{definition}
\newtheorem{dfn}[thm]{Definition} 
\newtheorem{ex}[thm]{Example}
\newtheorem{rmk}[thm]{Remark}

\numberwithin{thm}{section}
\numberwithin{equation}{section}

\newcommand{\Proof}{\noindent {\it Proof}.\ \ }
\newcommand{\Hom}{\operatorname{Hom}}
\newcommand{\Ext}{\operatorname{Ext}}
\newcommand{\Spec}{\operatorname{Spec}}

\newcommand{\ob}{\operatorname{ob}}
\newcommand{\Hilb}{\operatorname{Hilb}}
\newcommand{\red}{\operatorname{red}}
\newcommand{\im}{\operatorname{im}}
\newcommand{\Bs}{\operatorname{Bs}}

\newcommand{\Pic}{\operatorname{Pic}}

\newcommand{\Gr}{\operatorname{Gr}}
\newcommand{\Bl}{\operatorname{Bl}}
\newcommand{\algeeq}{\stackrel{alg.}{\; \sim \;}}

\newcommand{\car}{\operatorname{char}}

\renewcommand{\labelenumi}{{\rm (\arabic{enumi})}}

\newcommand{\mapright}[1]{%
\smash{\mathop{%
\hbox to 1cm{\rightarrowfill}}\limits^{#1}}}
\newcommand{\mapleft}[1]{%
\smash{\mathop{%
\hbox to 1cm{\leftarrowfill}}\limits^{#1}}}
\newcommand{\mapdown}[1]{\Big\downarrow
\llap{$\vcenter{\hbox{$\scriptstyle#1\,$}}$ }}
\newcommand{\mapup}[1]{\Big\uparrow
\rlap{$\vcenter{\hbox{$\scriptstyle#1$}}$ }}

\title{Obstructions to deforming curves on a $3$-fold, II: \\
{\Large Deformations of degenerate curves on 
a del Pezzo $3$-fold}}
\author{Hirokazu N{\sc asu}\thanks{
Supported in part by the 21-st century COE program
``Formation of an International Center of Excellence
in the Frontier of Mathematics
and Fostering of Researchers in Future Generations''.}}
\date{}

\begin{document}
\maketitle

\begin{abstract}
We study the Hilbert scheme $\Hilb^{sc} V$ of
smooth connected curves on a smooth del Pezzo $3$-fold $V$.
We prove that every degenerate curve $C$, i.e. 
every curve contained in a smooth hyperplane section $S$ of $V$,
does not deform to a non-degenerate curve 
if the following two conditions are satisfied:
(i) $\chi(V,\mathcal I_C(S))\ge 1$ and
(ii) for every line $\ell$ on $S$ such that $\ell \cap C = \emptyset$,
the normal bundle $N_{\ell/V}$ is trivial
(i.e.~$N_{\ell/V} \simeq {\mathcal O_{\mathbb P^1}}^{\oplus 2}$).
As a consequence, we prove an analogue 
(for $\Hilb^{sc} V$) of a conjecture of J.~O.~Kleppe 
which is concerned with non-reduced components of the Hilbert scheme
$\Hilb^{sc} \mathbb P^3$ of curves in 
the projective $3$-space $\mathbb P^3$.
\end{abstract}

\section{Introduction}\label{introduction}
This paper is a sequel to a joint work \cite{Mukai-Nasu} 
with Shigeru Mukai. In \cite{Mukai-Nasu} the embedded deformations of
smooth curves $C$ on a 
smooth projective $3$-fold $V$ have been studied under the presence of 
a smooth surface $S$ such that $C \subset S \subset V$,
especially when $V$ is a uniruled $3$-fold.
In this paper, the same subject is studied in detail especially 
when $V$ is a del Pezzo $3$-fold.

It is known that even if the deformations of 
$C$ in $S$ and the deformations of $S$ in $V$ behave well,
those of $C$ in $V$ behave badly in general.
For example, even if $\Hilb V$ and $\Hilb S$ are nonsingular 
of expected dimension $\chi(N_{S/V})$ and $\chi(N_{C/S})$ 
at $[S]$ and $[C]$ respectively, 
there can be a generically non-reduced component of 
$\Hilb V$ passing through $[C]$ 
(cf.~Mumford's example in \cite{Mumford}).
Such non-reduced components of
the Hilbert scheme $\Hilb^{sc} V$ of smooth connected curves on $V$
have been constructed for many uniruled $3$-folds $V$ in \cite{Mukai-Nasu}.
The non-reducedness is originated from the non-surjectivity of 
the restriction map
\begin{equation}\label{map:restriction}
H^0(S,N_{S/V}) \overset{|_C}\longrightarrow H^0(C,N_{S/V}\big{\vert}_C).
\end{equation}
We say that $C$ is {\em stably degenerate}
if every (small) deformation of $C$ in $V$
is contained in a divisor $S'$ of $V$
which is algebraically equivalent to $S$.
If \eqref{map:restriction} is surjective, 
then $C$ is stably degenerate 
(cf.~Proposition~\ref{prop:principle}).
If it is not surjective, then
there can be a first order deformation $\tilde C$ of $C$ in $V$
which is not contained in any first order deformation $\tilde S$ of $S$.
In this paper, we consider the following problem raised by Mukai:
\begin{prob}\label{prob:naive}
Suppose that \eqref{map:restriction} 
is not surjective and $\chi(V,\mathcal I_C(S))>0$.
Then (1) Is $C$ stably degenerate? (2) Is $\Hilb^{sc} V$ singular at $[C]$?
\end{prob}
Here $\mathcal I_C$ denotes the ideal sheaf of $C$ in $V$ and
$\mathcal I_C(S):=\mathcal I_C \otimes \mathcal O_V(S)$.
J.~O.~Kleppe \cite{Kleppe85} and Ph.~Ellia \cite{Ellia} 
considered Problem~\ref{prob:naive} for the case where
$V$ is the projective $3$-space $\mathbb P^3$,
$S$ is a smooth cubic surface in $\mathbb P^3$
and $C$ is a smooth connected curve on $S$.
Kleppe gave a conjecture (cf.~Conjectures~\ref{conj:Kleppe-Ellia}), 
which can be reformulated as follows:
\begin{conj}
\label{conj:Kleppe-Ellia-Nasu}
Let $C \subset S \subset \mathbb P^3$ be as above and assume that
$\chi(\mathbb P^3,\mathcal I_C(3))\ge 1$. 
If $C$ is linearly normal, 
then every (small) deformation $C'$ of $C$ in $\mathbb P^3$
is contained in a cubic surface $S' \subset \mathbb P^3$,
i.e.~$C$ is stably degenerate.
\end{conj}
As a testing ground of his conjecture,
we consider Problem~\ref{prob:naive} for the case where 
$V$ is a smooth del Pezzo $3$-fold (cf.~\S\ref{threefold}),
$S$ is a smooth polarization of $V$, i.e., a smooth member 
of the half anti-canonical system $|-\frac 12 K_V|$ and 
$C$ is a smooth connected curve on $S$.
The following theorem is an analogue of Kleppe's conjecture.

\begin{thm}\label{thm:main}
Let $C \subset S \subset V$ be as above and
assume that $\chi(V,\mathcal I_C(S))\ge 1$.
If every line $\ell$ on $S$ such that $C \cap \ell=\emptyset$
is a {\em good line} on $V$ 
(i.e.~the normal bundle $N_{\ell/V}$ of $\ell$ in $V$ is trivial), 
then:
\begin{enumerate}
 \item $C$ is stably degenerate, and
 \item $\Hilb^{sc} V$ is nonsingular at $[C]$
       if and only if $H^1(V,\mathcal I_C(S))= 0$.
\end{enumerate}
\end{thm}
If $\chi(V,\mathcal I_C(S))< 1$, then it follows from a dimension count
that $C$ is not stably degenerate 
(Proposition~\ref{prop:numerically non-degenerate}).
If some $\ell$ is a {\em bad line} on $V$
(i.e.~$N_{\ell/V}\not\simeq {\mathcal O_{\mathbb P^1}}^{\oplus 2}$)
then $C$ is not necessarily stably degenerate
(Proposition~\ref{prop:non-stably degenerate}).
As a corollary to Theorem~\ref{thm:main},
we give a sufficient condition for a maximal family $W$ of 
degenerate curves on $V$ to become 
an irreducible component of the Hilbert scheme 
$\Hilb^{sc} V$ and determine whether $\Hilb^{sc} V$ 
is generically non-reduced along $W$ or not
(Theorem~\ref{thm:generic smoothness}).

One of the main tools used in this paper 
is the infinitesimal analysis of the Hilbert scheme
developed in \cite{Mukai-Nasu}. As is well known, 
every infinitesimal deformation $\tilde C$ of $C$ in $V$
of first order (i.e.~over $\Spec k[t]/(t^2)$)
determines a global section $\alpha \in H^0(N_{C/V})$
and a cohomology class $\ob(\alpha) \in H^1(N_{C/V})$
called the {\em obstruction}.
Then $\tilde C$ lifts to a deformation over $\Spec k[t]/(t^3)$
if and only if $\ob(\alpha)=0$ (cf.~\S\ref{infinitesimal}).
Let $\pi_S: N_{C/V} \rightarrow N_{S/V}\big{\vert}_C$ be
the natural projection.
In \cite{Mukai-Nasu} Mukai and Nasu studied 
the {\em exterior component} of $\alpha$ and $\ob(\alpha)$,
i.e., the images of $\alpha$ and $\ob(\alpha)$
by the induced maps
$H^i(\pi_S): H^i(N_{C/V}) \rightarrow H^i(N_{S/V}\big{\vert}_C)$
($i=0,1$), respectively.
They proved that if there exists a curve $E$ on $S$ such that
$(E^2)_S<0$ (e.g. $(-1)$-$\mathbb P^1$ on $S$) and
the exterior component of $\alpha$ lifts to 
a global section $v \in H^0(N_{S/V}(E)) \setminus H^0(N_{S/V})$, 
then the exterior component of $\ob(\alpha)$ is nonzero 
provided that certain additional conditions on $E$, $C$ and $v$ hold
(see \cite[Theorem~1.6]{Mukai-Nasu}).
Such a rational section $v$ of $N_{S/V}$ 
admitting a pole along $E$
is called an {\em infinitesimal deformation with a pole}.
In \S\ref{with pole} we see that an infinitesimal deformation
with a pole along $E$ induces an obstructed infinitesimal deformation
of the open surface $S^{\circ}:=S\setminus E$ in 
the open $3$-fold $V^{\circ}:=V\setminus E$ 
(Theorem~\ref{thm:open surface}).
By using this fact, we prove Theorem~\ref{thm:main}
in \S\ref{obstruction}. 
In \S\ref{example} we give some examples of
generically non-reduced components of the Hilbert scheme
of curves on a del Pezzo $3$-fold as an application.

\paragraph{Acknowledgements}
I should like to express my sincere gratitude
to Professor Shigeru Mukai.
He showed me the example of 
non-reduced components of the Hilbert scheme
of canonical curves given in \S\ref{canonical}
as a simplification of Mumford's example of a non-reduced component 
of $\Hilb^{sc} \mathbb P^3$.
This led me to research the topic of this paper.
Throughout the research,
he made many suggestions which are useful 
for obtaining and improving the proofs.
In particular, according to his suggestion, 
I studied the embedded infinitesimal deformations of
a non-projective surface and organize \S\ref{with pole}
to improve the crucial part of the proof of Proposition~\ref{prop:core}.
I am grateful to Professor Jan Oddvar Kleppe for giving me
useful comments on Hilbert-flag schemes
and finding a gap in the original proof of Lemma~\ref{lem:nonspecial}.

\paragraph{Notation and Conventions}
We work over an algebraically closed field $k$ of characteristic $0$.
Let $V$ be a scheme over $k$ and let $X$ be a closed subscheme of $V$.
Then $\mathcal I_X$ denotes the ideal sheaf of $X$ in $V$ and
$N_{X/V}$ denotes the normal sheaf
$(\mathcal I_X /{\mathcal I_X}^2)^{\vee}$ of $X$ in $V$.
For a sheaf $\mathcal F$ on $V$,
we denote the restriction map 
$H^i(V,\mathcal F) \rightarrow H^i(X,\mathcal F\big{\vert}_X)$ 
by $\big{\vert}_X$.
We denote the Euler-Poincar\'e characteristic of $\mathcal F$ 
by $\chi(V,\mathcal F)$ or $\chi(\mathcal F)$.
$\Hilb^{sc} V$ denotes the open subscheme of the Hilbert scheme $\Hilb V$ 
whose point corresponds to a smooth connected curve on $V$.

\section{Preliminaries}\label{preliminaries}

\subsection{Del Pezzo surfaces}\label{surface}

A {\em del Pezzo surface} is a
smooth surface $S$ whose anti-canonical divisor $-K_S$ is ample.
Every del Pezzo surface is isomorphic to
$\mathbb P^2$ blown up at fewer than $9$ points or 
$\mathbb P^1 \times \mathbb P^1$.
We denote the blow-up of $\mathbb P^2$ at $(9-n)$-points by $S_n$.
A curve $\ell \simeq \mathbb P^1$ on $S_n$ 
is called a {\em line}\footnote{There exists no line on 
$\mathbb P^2$ and on $\mathbb P^1 \times \mathbb P^1$.}
if $\ell \cdot (-K_S)=1$.
Every $(-1)$-$\mathbb P^1$ on $S_n$ is a line 
and every line on $S_n$ is a $(-1)$-$\mathbb P^1$.
A curve $q$ on $S_n$ is called a {\em conic} if
$q \cdot (-K_S)=2$ and $q^2=0$.

\begin{lem}\label{lem:vanishing1}
 Let $D$ be a divisor on a del Pezzo surface $S$.
 If $D$ is nef and $\chi(-D)\ge 0$, then $H^1(S,-D)=0$.
\end{lem}
\Proof
If $D^2>0$ then the assertion follows the Kawamata-Viehweg vanishing.
Since $D$ is a nef divisor on a del Pezzo surface, we have $D^2 \ge 0$.
Now we assume that $D^2=0$.
If $S=S_n$, then $D$ is linearly equivalent to 
a multiple $mq$ ($m\ge 0$) of a conic $q$ on $S$.
By the Riemann-Roch theorem, we have
\begin{align*}
\chi(-D)
& =\dfrac12 (-mq)\cdot (-mq -K_S)+\chi(\mathcal O_S) \\
& =-m+1.
\end{align*}
Thus we have $m=0$ or $1$ by assumption. This implies 
that $H^1(-mq)=0$.
If $S=\mathbb P^1\times \mathbb P^1$, then
$D$ is of bidegree $(m,0)$ or $(0,m)$ with $m\ge 0$.
Again by the Riemann-Roch theorem, we have $\chi(-D)=-m+1\ge 0$.
Thus $H^1(\mathcal O_{\mathbb P^1\times \mathbb P^1}(-D))=0$.
\qed

\begin{lem}\label{lem:base}
 Let $D$ be an effective divisor on a del Pezzo surface $S$.
 Then the lines $\ell$ such that $D\cdot \ell < 0$ are mutually disjoint.
 The fixed part\footnote{%
 the base locus of dimension one}
 $\Bs |D|$ of the linear system $|D|$ on $S$
 is equal to
 $$
 - \sum_{D\cdot \ell < 0} (D\cdot \ell) \ell.
 $$ 
\end{lem}
\Proof
We prove the two assertions at the same time.
It is clear that any line $\ell$ satisfying $D\cdot \ell<0$ is contained 
in $\Bs |D|$. On the other hand, except for lines on $S$
every irreducible curve $C$ on $S$ can move on $S$ 
by the linearly equivalence since $\chi(C)\ge 2$ and $H^2(C)=0$.
Hence $|D|$ is decomposed into the sum
$$
|D| = |D'| + \sum_{i=1}^k m_i \ell_i, 
$$
of a linear system $|D'|$ on $S$ such that $\Bs|D'|=\emptyset$
and some lines $\ell_i$ on $S$
with coefficients $m_i \in \mathbb Z_{>0}$ ($1 \le i \le k$).
If $\ell_i \cap \ell_j \ne \emptyset$ for some $i\ne j$,
then $\ell_i+\ell_j$ is a (reducible) conic on $S$ and
can move on $S$ by $\chi(\ell_i+\ell_j)=2$.
Thus $\ell_i$'s are mutually disjoint.
Now we prove that $D \cdot \ell_i <0$ for any $i$.
Since $m_i=(D'-D)\cdot \ell_i>0$,
it suffices to show that $D'\cdot \ell_i=0$.
Since $D'$ is nef, we have $(D')^2\ge 0$.
Since $-K_S$ is ample, so is $D'-K_S$.
Hence we have $H^1(D')=H^1((D'-K_S)+K_S)=0$ by 
the Kodaira vanishing.
If $D'\cdot \ell_i\ge 1$, then it follows from the exact sequence
$$
0 \longrightarrow \mathcal O_S(D')
\longrightarrow \mathcal O_S(D'+\ell_i)
\longrightarrow \mathcal O_S(D'+\ell_i)\big{\vert}_{\ell_i}
\longrightarrow 0
$$
that $h^0(D'+\ell_i)> h^0(D')$. Thus we have $D'\cdot \ell_i =0$.
\qed

\begin{lem}\label{lem:freeness}
 Let $E$ be a disjoint union of $m$ lines ($m \ge 0$) 
 on a del Pezzo surface $S$ and
 let $\varepsilon :S \rightarrow F$ be the blow-down of $E$ from $S$.
 If a divisor $D$ on $F$ satisfies $h^0(F,D)\ge m$,
 then we have the following:
 \begin{enumerate}
  \item $h^0(S,\varepsilon^*D-E)=h^0(F,D)-m$, and
  \item If $H^1(S,\varepsilon^*D)=0$, 
	then $H^1(S,\varepsilon^* D-E)=0$.
 \end{enumerate}
\end{lem}
\Proof
(1) Let $\ell_i$ $(1\le i\le m)$ be the disjoint lines on $S$ and
let $E:=\sum_{i=1}^m \ell_i$.
We put $D_j:=\varepsilon^* D-\sum_{1\le i\le j} \ell_i$.
Since the image of $\ell_i$ on $F$ is a point, 
we have $h^0(D_j) \ge h^0(D) -j$ for every $1 \le j \le m$.
Moreover since $D_{j-1} \cdot \ell_j=0$, Lemma~\ref{lem:base}
shows that $\ell_j$ is not contained in $\Bs|D_{j-1}|$.
Hence $\dim |D_j|$ decreases one by one as $j$ increases.
Therefore we have $h^0(\varepsilon^*D-E)=h^0(D_m)=h^0(D)-m$.

(2) An exact sequence
$0 \rightarrow \mathcal O_S(\varepsilon^* D-E)
\rightarrow \mathcal O_S(\varepsilon^* D)
\rightarrow \mathcal O_E
\rightarrow 0$ on $S$ induces an exact sequence
$$
H^0(S,\varepsilon^*D)\overset{\rho}\longrightarrow H^0(E,\mathcal O_E)
\longrightarrow H^1(S,\varepsilon^*D-E) 
\longrightarrow H^1(S,\varepsilon^*D)
$$
of cohomology groups. Then $\rho$ is surjective by (1)
and $H^1(S,\varepsilon^* D)=0$ by assumption.
Hence we have $H^1(S,\varepsilon^* D-E)=0$.
\qed

\medskip

Let $C$ be a smooth connected curve on a del Pezzo surface $S$.
We consider the restriction to $C$ of the anti-canonical 
linear system $|-K_S|$ on $S$. The restriction map
$H^0(-K_S) \rightarrow H^0(-K_S\big{\vert}_C)$ 
is not surjective in general. 
Let $\ell_i$ ($1 \le i \le m$) be the lines on $S$
disjoint to $C$.
Let us define an effective divisor $E$ on $S$ by the sum
$$
E := \sum_{i=1}^m \ell_i
$$
and we put $E:=0$ if there exists no such $\ell_i$.
If $C$ is neither a line nor a conic, 
then $\ell_i$'s are mutually disjoint:
indeed if $\ell_i \cap \ell_j \ne \emptyset$ for some $i \ne j$,
then $q:=\ell_i + \ell_j$ is a conic on $S$ 
and hence $C$ intersects $q$ by $C \cdot q>0$.

\begin{prop}\label{prop:anti-canonical}
 Assume that $C$ is not rational and $\chi(-K_S-C)\ge 0$.
 Then we have $H^1(S,-K_S+E-C)=0$ and the restriction map
  \begin{equation}\label{map:anti-canonical}
   H^0(S,-K_S+E) \overset{|_C}{\longrightarrow}
   H^0(C,-K_S\big{\vert}_C) 
  \end{equation}
 is surjective. If $C$ is not elliptic either,
 then the map \eqref{map:anti-canonical} is an isomorphism.
\end{prop}
\Proof
It suffices to show that $H^1(-K_S+E-C)=0$ by 
the exact sequence
\begin{equation}\label{ses:restriction to C}
0 \longrightarrow \mathcal O_S(-K_S+E-C)
\longrightarrow \mathcal O_S(-K_S+E)
\longrightarrow \mathcal O_S(-K_S)\big{\vert}_C
\longrightarrow 0. 
\end{equation}

\noindent
{\bf Claim.} \qquad
Put $D_1:=C+K_S-E$. Then $D_1$ is nef.

\medskip

Since $S$ is regular (i.e.~$H^1(K_S)=0$),
the restriction map 
$\big{\vert}_C: H^0(C+K_S)\rightarrow H^0(K_C)$ 
is surjective. Since $C \not\simeq \mathbb P^1$, 
the linear system $|C+K_S|$ on $S$ is non-empty. 
Let $l$ be a line on $S$. Since $C$ is not a line, 
we have $C\cdot \ell \ge 0$ and hence $(C+K_S)\cdot \ell \ge -1$.
By Lemma~\ref{lem:base}, $\ell$ is contained in $\Bs |C+K_S|$
if and only if $C \cap \ell=\emptyset$.
Thus we have $E=\Bs |C+K_S|$ and $|D_1|$ does not have base components.
In particular, $D_1$ is nef.

\medskip

It follows from the exact sequence 
\begin{equation}\label{ses:restriction to E}
0 \longrightarrow \mathcal O_S(-K_S-C)
\longrightarrow \mathcal O_S(-K_S+E-C) \longrightarrow 
\underbrace{\mathcal O_S(-K_S+E)\big{\vert}_E}_{\simeq \, \mathcal O_E}
\longrightarrow 0 
\end{equation}
that $\chi(-D_1)=\chi(-K_S-C)+\chi(\mathcal O_E)\ge 0$.
Hence we have $H^1(-D_1)=0$ by Lemma~\ref{lem:vanishing1}.

Now we assume that $C$ is not elliptic.
Then $K_C \not\sim 0$ and hence $C+K_S \not\sim E$ 
by adjunction. Thus $D_1 \not\sim 0$ and $H^0(-D_1)=0$.
Therefore \eqref{map:anti-canonical} is injective.
\qed

\begin{lem}\label{lem:injective2}
 If $C$ is not rational nor elliptic and $\chi(-K_S-C)\ge 0$, 
 then the map 
 $$
 H^1(S,-K_S+3E) \overset{|_C}\longrightarrow 
 H^1(C,-K_S\big{\vert}_C)
 $$
 induced by \eqref{ses:restriction to C}$\otimes \mathcal O_S(2E)$
 is injective.
\end{lem}
\Proof
It suffices to show that $H^1(-K_S+3E-C)=0$.
Let $\varepsilon: S \rightarrow F$ be the blow-down of $E$ from $S$.
Then there exists a divisor $D_2$ on $F$ such that
$\varepsilon^*D_2 \sim C+2K_S-2E$.
By the Serre duality, it suffices to show that
$H^1(\varepsilon^*D_2-E)=0$.

\medskip

\noindent
{\bf Claim.} \qquad
$H^i(S,\varepsilon^*D_2)=0$ for $i =1,2$.

\medskip

By \eqref{ses:restriction to E}$\otimes \mathcal O_S(E)$,
there exists an exact sequence 
$$
H^1(S,-K_S+E-C) \longrightarrow H^1(S,-K_S+2E-C) \longrightarrow 
H^1(E,(-K_S+2E)\big{\vert}_E).
$$
Since $H^1((-K_S+2E)\big{\vert}_E)\simeq H^1(\mathcal O_E(E))=0$ and
$H^1(-K_S+E-C)=0$ by Proposition~\ref{prop:anti-canonical},
we have $H^1(-K_S+2E-C)=0$.
By the Serre duality, we have $H^1(\varepsilon^*D_2)=0$.
Similarly by the Serre duality, we have
$H^2(\varepsilon^*D_2)\simeq H^0(K_S-\varepsilon^*D_2)^{\vee}$.
Since $C$ is not rational nor elliptic, we have
$(K_S-\varepsilon^*D_2)\cdot C=(-K_S-C)\cdot C=-\deg K_C<0$.
Hence we have $H^2(\varepsilon^*D_2)=0$ because $C$ is nef.
Thus the claim has been proved.

\medskip

By this claim, we have
$h^0(F,D_2)=h^0(S,\varepsilon^*D_2)= \chi(S,\varepsilon^*D_2)$.
Then an easy calculation shows that
$\chi(\varepsilon^*D_2)=\chi(-K_S-C)+\chi(\mathcal O_E)$.
Since $\chi(-K_S-C)\ge 0$,
we have $h^0(F,D_2)=\chi(S,\varepsilon^*D_2)\ge m$, 
where $m$ is the number of components of $E$.
Since $H^1(\varepsilon^*D_2)=0$,
Lemma~\ref{lem:freeness}~(2) shows that $H^1(\varepsilon^*D_2-E)=0$.
\qed

\medskip

Let $S$ be a smooth projective surface and let $L$ be a line bundle on $S$.

\begin{lem}\label{lem:injective}
 Let $E$ be a disjoint union of 
 irreducible curves $E_i$ ($i=1,\ldots,m$) on $S$ such that $E_i^2<0$
 and let $\iota: S^{\circ}:=S \setminus E \hookrightarrow S$ be the open 
 immersion. If $\deg(L\big{\vert}_{E_i}) \le 0$ for every $i$, 
 then the map
 $$
 H^1(S,L) \rightarrow H^1(S^{\circ},L\big{\vert}_{S^{\circ}})
 $$
 induced by the sheaf inclusion 
 $L \hookrightarrow L\otimes \iota_*\mathcal O_{S^{\circ}}$ is injective.
\end{lem}
The proof is similar to that of {\cite[Lemma~2.3]{Mukai-Nasu}}
and we omit it here.
Lemma~\ref{lem:injective} allows us to identify 
$H^1(S,L(nE))$ ($n \ge 0$) with their images in 
$H^1(S^{\circ},L\big{\vert}_{S^{\circ}})$.
As a result, under the identification we obtain a natural filtration
$$
H^1(S,L) \subset
H^1(S,L(E)) \subset
H^1(S,L(2E)) \subset \dots
\subset H^1(S^{\circ},L\big{\vert}_{S^{\circ}})
$$
on $H^1(S^{\circ},L\big{\vert}_{S^{\circ}})$.

\subsection{Del Pezzo threefolds}\label{threefold}

A {\em del Pezzo threefold} is a pair $(V,H)$ consisting of 
a (smooth) irreducible projective variety $V$ of dimension $3$
and an ample Cartier divisor $H$ on $V$ such that $-K_V=2H$.
Here $H$ is called the {\em polarization} of $V$ and sometimes omitted.
The self-intersection number $n:=H^3$ is called the
{\em degree} of $V$. It is known that
the linear system $|H|$ on $V$ determines a double cover
$\varphi_{|H|}:V \rightarrow \mathbb P^3$ if $n=2$,
and an embedding $\varphi_{|H|}:V \hookrightarrow \mathbb P^{n+1}$
if $n\ge 3$. 
If $S$ is a smooth member of $|H|$, then the pair $(S,H\big{\vert}_S)$
is a del Pezzo surface of degree $n$.
Every smooth del Pezzo $3$-fold is one of 
$V_n$ ($1 \le n \le 8$) or $V_6'$ in Table~\ref{table:delpezzo 3-fold},
\begin{table}[h]
\begin{minipage}{16cm}
\caption{Del Pezzo $3$-folds}
\label{table:delpezzo 3-fold}
\begin{center}
\begin{tabular}{|l|c|c|l|}
\hline
del Pezzo $3$-folds & $n$ & $\rho$ &  \\
\hline
 $V_1 = (6) \subset \mathbb P(3,2,1,1,1)$ & $1$ & $1$ & 
 a weighted hypersurface of degree $6$ \\
 \hline
 $V_2 = (4) \subset \mathbb P(2,1,1,1,1)$ & $2$ & $1$ & 
 a weighted hypersurface of degree $4$
 \footnote{Another realization of $V_2$ is a double cover of 
 $\mathbb P^3$ branched along a quartic surface.}
 \\
 \hline
 $V_3 = (3) \subset \mathbb P^4$ & $3$ & $1$ & a cubic hypersurface \\
 \hline
 $V_4 = (2) \cap (2) \subset \mathbb P^5$ & $4$ & $1$ & 
 a complete intersection of two quadrics \\
 \hline
 $V_5 = [\Gr(2,5) \overset{\mbox{\scriptsize
 Pl{\"u}cker}}{\hookrightarrow} \mathbb P^9] 
 \cap \mathbb L^{(6)}$ & $5$ & $1$ & 
 a linear section of Grassmannian\\
 \hline
 $V_6= [\mathbb P^1 \times \mathbb P^1 \times \mathbb P^1 
 \overset{\rm Segre}{\hookrightarrow} \mathbb P^7]$  & $6$ & $3$ & \\
 \hline 
 $V_6' = [\mathbb P^2 \times \mathbb P^2 \overset{\rm
 Segre}{\hookrightarrow} \mathbb P^8] \cap \mathbb L^{(7)}$ & $6$ & $2$ & \\
 \hline
 $V_7 = \Bl_{\rm pt} \mathbb P^3 \subset \mathbb P^8$ & $7$ & $2$ & the blow-up
 of $\mathbb P^3$ at a point
 \footnote{$V_7$ is realized as the projection of 
 $V_8 \subset \mathbb P^9$ from one of its point.}
 \\
 \hline
 $V_8 = \mathbb P^3 \overset{\rm Veronese}{\hookrightarrow} \mathbb P^9$ & $8$
 & $1$ & the Veronese image of $\mathbb P^3$ \\
 \hline
\end{tabular}
\end{center}  
\end{minipage}
\end{table}
in which $\mathbb L^{(i)}$ denotes a linear subspace of dimension $i$,
and $n$ and $\rho$ respectively denote the degree and 
the Picard number of $V_n$ (and of $V_6'$)
(cf.~\cite{Fujita80},\cite{Fujita81},\cite{Iskovskih77}).
It is known that a smooth $3$-fold $V \subset \mathbb P^{n+1}$ 
($n \ge 3$) is a del Pezzo $3$-fold of degree $n$
if a linear section
$[V \subset \mathbb P^{n+1}] \cap H_1 \cap H_2$
with two general hyperplanes $H_1,H_2 \subset \mathbb P^{n+1}$
is an elliptic normal curve in $\mathbb P^{n-1}$.

We briefly review the basics of the Hilbert scheme of lines on a
del Pezzo $3$-fold. We refer to 
Iskovskih (\cite{Iskovskih77},\cite{Iskovskih79})
for the details.
Let $(V,H)$ be a smooth del Pezzo $3$-fold of degree $n$.
By a {\em line} on $(V,H)$, we mean a reduced irreducible curve 
$\ell$ on $V$ such that $(\ell\cdot H)_V=1$ and $\ell \simeq \mathbb P^1$.
If $n \le 7$ then $V$ contains a line $\ell$.
Then there are only the following possibilities for the normal 
bundle $N_{\ell/V}$ of $\ell$ in $V$:
 $$
 \begin{array}{cll}
  \mbox{(0,0):} & 
   N_{\ell/V} \simeq {\mathcal O_{\mathbb P^1}}^{\oplus 2}
   & (\mbox{i.e. trivial}), \\
  \mbox{(1,-1):} &
   N_{\ell/V} \simeq \mathcal O_{\mathbb P^1}(-1)
   \oplus \mathcal O_{\mathbb P^1}(1), & \\
  \mbox{(2,-2):} &
   N_{\ell/V} \simeq \mathcal O_{\mathbb P^1}(-2)
   \oplus \mathcal O_{\mathbb P^1}(2) & (\mbox{only if $n=1$ or $2$}), \\  
  \mbox{(3,-3):} &
   N_{\ell/V} \simeq \mathcal O_{\mathbb P^1}(-3)
   \oplus \mathcal O_{\mathbb P^1}(3) & (\mbox{only if $n=1$}). \\
 \end{array}
 $$
In this paper,
$\ell$ is called a {\em good line} if $N_{\ell/V}$ is trivial,
and called a {\em bad line} otherwise.
If $n \ge 3$, then every line on $V$ is of type $(0,0)$ or $(1,-1)$.
The Hilbert scheme $\Gamma$ of lines on $V$ is 
called the {\em Fano surface} of $V$, and in fact
every irreducible (non-embedded)
component of $\Gamma$ is of dimension two.
Let $\Gamma_i \subset \Gamma$ be an irreducible component
and let $S_i$ be the universal family of lines on $V$ over $\Gamma_i$.
Then there exists a natural diagram.
$$
\begin{array}{ccc}
 S_i & \mapright{p} & V \\
 \mapdown{\pi} & & \\
 \Gamma_i. &&
\end{array}
$$
By \cite[Chap.III, Proposition~1.3 (iv)]{Iskovskih79}, 
if $n \ge 3$ then we have either
\begin{enumerate}
 \renewcommand{\labelenumi}{{\rm (\alph{enumi})}}
 \item $p$ is surjective; in this case a general line in
       $\Gamma_i$ is a good line; or
 \item $p(S_i)\simeq \mathbb P^2$ is a plane on $V \subset \mathbb P^{n+1}$;
       in this case every line in $\Gamma_i$ is a bad line.
\end{enumerate}
We have either (a) or (b) also when $n\le 2$.
(See the proof\footnote{%
In the proof, the assumption that $\car k=0$ is used.}
in \cite{Iskovskih79}, 
which works for $n \le 2$.)
If $n \ne 7$ then every irreducible component of $\Gamma$ is of
type $(a)$. If $n=7$ then $\Gamma$ consists of two 
irreducible components
$\Gamma_i \simeq \mathbb P^2 (i=0,1)$, one of which is of type (a),
while the other is of type (b). Consequently, we have
\begin{lem}[Iskovskih]\label{lem:good line}
Every smooth del Pezzo $3$-fold of degree $n \ne 8$ contains a good line.
\end{lem}

\begin{lem}\label{lem:general section}
 Let $(V,H)$ be a smooth del Pezzo $3$-fold of degree $n$ and let 
 $S$ be a general member of $|H|$.
 If $n\ne 7$ then $S$ does not contain a bad line.
 If $n = 7$ then $S$ contains three lines, one of which is bad,
 while the others are good.
\end{lem}
\Proof
There exists no line on $V_8$. 
If $n\ne 7$, then the locus 
$\mathfrak B$ of bad lines in the Fano surface $\Gamma$ is of dimension one.
Let $p_i$ denote the projection of
$$
\left\{(\ell,S) \bigm| \ell \subset S \right\} \subset \Gamma \times |H|
$$
to the $i$-th factor. Since the fiber of $p_1$ is of dimension $n-1$,
$p_2(p_1^{-1}(\mathfrak B))$ is a proper closed subset
of $|H|\simeq \mathbb P^{n+1}$.
Hence every line on a general member $S$ of $|H|$ is a good line.

Suppose that $V=V_7$, i.e., the blow-up of $\mathbb P^3$ at a point.
Then $S$ is a del Pezzo surface $S_7$, i.e.,
a blow-up of $\mathbb P^2$ at two distinct points. 
Hence there are three lines
(i.e.~three $(-1)$-$\mathbb P^1$'s)
$\ell_0,\ell_1,\ell_2$ on $S$ forming the configuration
in Figure~\ref{fig:lines on S_7}:

\begin{figure}[h]
\begin{center}
\begin{picture}(50,50)(0,0)
 \put(10,5){\line(0,1){40}}
 \put(40,5){\line(0,1){40}}
 \put(0,10){\line(1,0){50}}
 \put(-10,30){$\ell_1$}
 \put(45,30){$\ell_2$}
 \put(20,-5){$\ell_0$}
\end{picture} 
 \caption{$(-1)$-$\mathbb P^1$'s on $S_7$}
\label{fig:lines on S_7}
\end{center}
\end{figure}
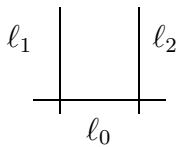
Here $\ell_0$ is distinguished by the fact that it intersects both of
the other lines.
Let $P$ be the exceptional divisor of the blow-up 
$V_7 \rightarrow \mathbb P^3$.
Then $P \simeq \mathbb P^2$ is a unique plane on $V_7$
and $\ell_0$ is the intersection of $S$ with $P$
(cf.~\cite[Chap II, \S1.4]{Iskovskih79}).
Since $N_{\ell_0/P} \simeq \mathcal O_{\mathbb P^1}(1)$,
$\ell_0$ is a bad line on $V_7$.
On the other hand, $\ell_1$ and $\ell_2$ are good lines on $V_7$
since $S$ is general.
\qed

\subsection{Infinitesimal deformations and  obstructions}
\label{infinitesimal}
Let $V$ be a smooth variety and let $X$ be a smooth closed subvariety of $V$.
An {\em (embedded) first order infinitesimal deformation} 
of $X$ in $V$ is a closed subscheme
$\tilde X \subset V\times \Spec k[t]/(t^2)$ 
which is flat over $\Spec k[t]/(t^2)$ and 
whose central fiber is $X$.
It is well known that 
there exists a one to one correspondence between
the group of homomorphisms $\alpha: \mathcal I_X \rightarrow \mathcal O_X$
and the first order infinitesimal deformations $\tilde X$ of $X$ in $V$.
In what follows, we identify $\tilde X$ with $\alpha$ and 
abuse the notation. The standard exact sequence
\begin{equation}\label{ses:standard sequence}
0 \longrightarrow \mathcal I_X
  \longrightarrow \mathcal O_V
  \longrightarrow \mathcal O_X
  \longrightarrow 0
\end{equation}
induces 
$\delta: \Hom(\mathcal I_X,\mathcal O_X)
\rightarrow \Ext^1(\mathcal I_X,\mathcal I_X)$
as a coboundary map.
Then $\alpha \in \Hom(\mathcal I_X,\mathcal O_X)$ (i.e.~$\tilde X$)
lifts to a deformation over $\Spec k[t]/(t^3)$
if and only if 
 $$
 \ob(\alpha):= \delta(\alpha) \cup \alpha 
 \in \Ext^1(\mathcal I_X,\mathcal O_X)
 $$
is zero,
where $\cup$ is the cup product map
$$
\Ext^1(\mathcal I_X,\mathcal I_X) \times
\Hom(\mathcal I_X,\mathcal O_X)
\overset{\cup}{\longrightarrow} \Ext^1(\mathcal I_X,\mathcal O_X).
$$
(We refer to \cite[Chap.~I \S2]{Kollar}. See also
\cite{Nasu}, \cite{Curtin}, \cite{Floystad} and \cite{Kleppe06}.)
Then $\ob(\alpha)$ is called the {\em obstruction} of 
$\alpha$ (i.e.~$\tilde X$).
Since both $X$ and $V$ are smooth, $\ob(\alpha)$ is contained
in $H^1(X,N_{X/V})\subset \Ext^1(\mathcal I_X,\mathcal O_X)$
(cf.~\cite[Chap.~I, Prop.~2.14]{Kollar}).
Since $\Hom(\mathcal I_X,\mathcal O_X) \simeq H^0(N_{X/V})$,
we regard $\alpha$ as a global section of $N_{X/V}$ from now on.

If $X$ is a hypersurface of $V$, i.e., of codimension one in $V$,
then $\ob(\alpha)$ becomes a simple cup product.
Let $\delta_1: H^0(X,N_{X/V}) \rightarrow H^1(V,\mathcal O_V)$
be the coboundary map of the exact sequence
$0 \rightarrow \mathcal O_V \rightarrow \mathcal O_V(X)
\rightarrow N_{X/V} \rightarrow 0$.
Let us define a map
\begin{equation}\label{map:d-map}
d_X: H^0(X,N_{X/V}) \longrightarrow H^1(X,\mathcal O_X) 
\end{equation}
by the composition of $\delta_1$ and the restriction map
$H^1(\mathcal O_V) \overset{\vert_X}{\longrightarrow} 
H^1(\mathcal O_X)$.\footnote{%
The map $d_X$ is equal to the map 
$d_{X,\mathcal O_V(X)}$ defined in \cite[\S2.1]{Mukai-Nasu}.
If $V$ is projective, then $d_X$ is the tangential map of
a natural morphism $X' \mapsto \mathcal O_X(X')$  from the
Hilbert scheme of divisor $X' \subset V$ to the Picard scheme
$\Pic X$.}
Then we have

\begin{lem}\label{lem:divisorial}
Let $X$ be a smooth hypersurface of $V$.
Then $\ob(\alpha)$ for $\alpha \in H^0(N_{X/V})$
is equal to the cup product $d_X(\alpha) \cup \alpha$, 
where $\cup$ is the cup product map
$$
H^1(X,\mathcal O_X) \times H^0(X,N_{X/V})
\overset{\cup}{\longrightarrow}
H^1(X,N_{X/V}).
$$ 
\end{lem}
\Proof
Since $\mathcal I_X \simeq \mathcal O_V(-X)$ is a line bundle on $V$,
we have
$\Ext^i(\mathcal I_X,\mathcal O_X) \simeq
H^i(N_{X/V})$ ($i=0,1$) and
$\Ext^1(\mathcal I_X,\mathcal I_X) \simeq
H^1(\mathcal O_V)$.
Hence the coboundary map $\delta$ appearing 
in the definition of $\ob(\alpha)$ is nothing but the 
coboundary map $\delta_1$ of
\eqref{ses:standard sequence}$\otimes \mathcal O_V(X)$.
Since $\alpha$ is a cohomology class on $X$,
the cup product map 
$H^1(\mathcal O_V) \rightarrow H^1(N_{X/V})$ with $\alpha$
factors through the restriction map $\big{\vert}_X$.
\qed

\medskip

We recall the definition of 
exterior component introduced in \cite{Mukai-Nasu}.
Let $X$ be a smooth closed subvariety of $V$ and 
let $Y$ be a smooth hypersurface of $V$ containing $X$.
Then the natural projection 
$\pi_Y: N_{X/V} \rightarrow N_{Y/V}\big{\vert}_X \simeq \mathcal O_X(Y)$
of normal bundles induces the maps
$H^i(\pi_Y): H^i(N_{X/V}) \rightarrow H^i(N_{Y/V}\big{\vert}_X)$,
where $i=0,1$, of their cohomology groups.
Let $\alpha$ be a global section of $N_{X/V}$.

\begin{dfn}\label{dfn:exterior}
$\pi_Y(\alpha)$ and $\ob_Y(\alpha)$ denote
the images of $\alpha$ and $\ob(\alpha)$
by the maps $H^0(\pi_Y)$ and $H^1(\pi_Y)$, respectively.
They are called the {\em exterior components} of $\alpha$ and
 $\ob(\alpha)$, respectively.
\end{dfn}

Roughly speaking, $\pi_Y(\alpha)$ is the projection of 
the normal vector $\alpha$ of $X$ in $V$
onto the normal directions to $Y$ in $V$.
Then $\ob_Y(\alpha)$ represents the obstruction to deforming 
$X$ into this directions.
We recall a basic fact on exterior components.

\begin{lem}[{\cite[Lemma~2.2]{Mukai-Nasu}}]\label{lem:reduction}
 Let $\pi_Y(\alpha)$ and $\ob_Y(\alpha)$ be 
 the exterior components of $\alpha$ and $\ob(\alpha)$, respectively.
 If there exists a global section $v$ of $N_{Y/V}$
 whose restriction $v\big{\vert}_X$ to $X$ 
 coincides with $\pi_Y(\alpha)$, then we have
 $$
 \ob_Y(\alpha) = \ob(v)\big{\vert}_X
 $$
 where $\ob(v)\big{\vert}_X \in H^1(N_{Y/V}\big{\vert}_X)$ 
 is the restriction of $\ob(v)\in H^1(N_{Y/V})$ to $X$.
\end{lem}
Lemma~\ref{lem:reduction} together with Lemma~\ref{lem:divisorial}
shows that $\ob_Y(\alpha)=d_Y(v)\big{\vert}_X \cup \pi_Y(\alpha)$,
where $d_Y$ is the map \eqref{map:d-map} for $Y$ and
$\cup$ is the cup product map
\begin{equation}\label{map:cup product}
 H^1(X,\mathcal O_X) \times H^0(X,N_{Y/V}\big{\vert}_X)
  \overset{\cup}\longrightarrow H^1(X,N_{Y/V}\big{\vert}_X). 
\end{equation}

Let $E$ be an effective divisor of $Y$ disjoint to $X$
(i.e.~$X \cap E =\emptyset$).
Let $Y^{\circ}$ and $V^{\circ}$ denote the two complements 
of $E$ in $Y$ and $V$, respectively.
Every rational section $v$ of $N_{Y/V}\simeq \mathcal O_Y(Y)$
having poles only along $E$ determines a global section 
$v^{\circ}$ of the normal sheaf $N_{S^{\circ}/V^{\circ}}$
of $Y^{\circ}$ in $V^{\circ}$
and hence obstruction 
$\ob(v^{\circ}) \in H^1(N_{S^{\circ}/V^{\circ}})$ to deforming
$S^{\circ}$ in $V^{\circ}$.
Let $\iota$ denote the open immersion of $Y^{\circ} \hookrightarrow Y$.
Then a natural homomorphism
$\iota_* N_{Y^{\circ}/V^{\circ}} \rightarrow N_{Y/V}\big{\vert}_X$
($=[\iota_* \mathcal O_{Y^{\circ}} \rightarrow \mathcal O_X]
\otimes N_{Y/V}$) of sheaves on $Y$ induces a map
$H^1(N_{S^{\circ}/V^{\circ}}) 
\overset{|_X}\longrightarrow H^1(N_{S/V}\big{\vert}_X)$.
Since $\ob(\alpha)$ is (and hence $\ob_Y(\alpha)$ is) determined
by a neighborhood of $X$, we have the following
variant of Lemma~\ref{lem:reduction}.

\begin{lem}\label{lem:reduction2}
 Let $\alpha$ be a global section of $N_{X/V}$.
 If there exists a rational section $v$ of $N_{S/V}$
 whose only poles are along $E$ and
 whose restriction to $X$ coincides with $\pi_Y(\alpha)$,
 then we have
 $$
 \ob_Y(\alpha)=\ob(v^{\circ}) \big{\vert}_X,
 $$
 where $\ob(v^{\circ})\big{\vert}_X$ is the image of
 $\ob(v^{\circ})$ by the map 
 $H^1(N_{S^{\circ}/V^{\circ}}) 
 \overset{|_X}\longrightarrow H^1(N_{S/V}\big{\vert}_X)$.
\end{lem}

\section{Infinitesimal deformations with a pole}\label{with pole}
Let $V$ be a smooth projective $3$-fold, $S$ a smooth surface in $V$, 
$E$ a smooth curve on $S$.
We put $V^{\circ}:=V\setminus E$ and $S^{\circ}:=S\setminus E$, 
the complemental open subvarieties.
In this section, we study the first order 
infinitesimal deformations of 
$S^{\circ}$ in $V^{\circ}$,
when the self-intersection number of $E$ on $S$ is negative.
We are interested in a rational section $v$ of $N_{S/V}$
having a pole only along $E$ and of order one,
that is, $v \in H^0(N_{S/V}(E))\setminus H^0(N_{S/V})$.
Let $\iota: S^{\circ} \hookrightarrow S$ be the open immersion.
Then $\iota_*\mathcal O_{S^{\circ}}$ contains
$\mathcal O_S(nE)$ as a subsheaf for any $n \ge 0$.
Hence the natural sheaf injection 
$N_{S/V}(nE) \hookrightarrow \iota_* N_{S^{\circ}/V^{\circ}}$
induces $H^0(S,N_{S/V}(nE)) \hookrightarrow 
H^0(S^{\circ},N_{S^{\circ}/V^{\circ}})$ for each $n$.
Therefore $v$ determines a first order infinitesimal deformation 
of $S^{\circ}$ in $V^{\circ}$.
The main theorem of this section is the following.

\begin{thm}\label{thm:open surface}
 Let $v$ be as above and assume that 
 $E^2<0$ and $\det N_{E/V}:=\bigwedge^2 N_{E/V}$ is trivial. 
 If the exact sequence
 \begin{equation}\label{ses:normal bundle of line}
   0 \longrightarrow N_{E/S} \longrightarrow 
   N_{E/V} \longrightarrow N_{S/V}\big{\vert}_E
   \longrightarrow 0
 \end{equation}
 does not split, then the first order infinitesimal deformation of 
 $S^{\circ} \subset V^{\circ}$ determined by $v$
 does not lift to a deformation over $\Spec k[t]/(t^3)$.
\end{thm}

Let $n$ be a non-negative integer. In what follows,
we identify $H^0(N_{S/V}(nE))$ with its image in 
$H^0(N_{S^{\circ}/V^{\circ}})$.
We shall prove that the obstruction
$\ob(v)$ is nonzero in $H^1(N_{S^{\circ}/V^{\circ}})$.
Let $d_{S^{\circ}}$ denote the map \eqref{map:d-map} for $X=S^{\circ}$.
Then by Lemma~\ref{lem:divisorial},
$\ob(v)$ is equal to the cup product 
$d_{S^{\circ}}(v) \cup v$,
where $\cup$ is the cup product map
$$
H^1(S^{\circ}, \mathcal O_{S^\circ}) \times 
H^0(S^{\circ},N_{S^{\circ}/V^{\circ}}) \overset{\cup}\longrightarrow 
H^1(S^{\circ},N_{S^{\circ}/V^{\circ}}).
$$
The inclusion
$\mathcal O_S(nE) \hookrightarrow \iota_* \mathcal O_{S^{\circ}}$
of sheaves induces a map $H^1(S,\mathcal O_S(nE)) \rightarrow
H^1(S^{\circ},\mathcal O_{S^{\circ}})$
of cohomology groups.
Suppose that $E^2 <0$.
Then this map is injective by Lemma~\ref{lem:injective}.
Hence we identify $H^1(\mathcal O_S(nE))$ with its image in
$H^1(S^{\circ},\mathcal O_{S^{\circ}})$.
Under this identification, there exists a natural filtration
$$
H^1(S,\mathcal O_S) \subset H^1(S,\mathcal O_S(E)) \subset
H^1(S,\mathcal O_S(2E)) \subset \cdots \subset 
H^1(S^{\circ},\mathcal O_{S^{\circ}})
$$
on $H^1(S^{\circ},\mathcal O_{S^{\circ}})$.
Suppose now that $\det N_{E/V}$ is trivial.
Then under similar identifications, there exists a natural filtration
$$
H^1(S,N_{S/V}(E)) \subset
H^1(S,N_{S/V}(2E)) \subset \cdots \subset 
H^1(S^{\circ},N_{S^{\circ}/V^{\circ}})
$$
on $H^1(S^{\circ},N_{S^{\circ}/V^{\circ}})$,
because we have $\deg N_{S/V}(nE)\big{\vert}_E
=\deg (\det N_{E/V}) + (n-1) E^2= (n-1)E^2 \le 0$
for $n \ge 1$.
Then it follows from \cite[Proposition~2.4~(1)]{Mukai-Nasu}
that the image of $d_{S^{\circ}}$ over $H^0(N_{S/V}(E))$
is contained in $H^1(\mathcal O_S(2E))$.
By the commutative diagram
$$
\begin{array}{ccccc}
 H^1(\mathcal O_{S^{\circ}}) & \times & H^0(N_{S^{\circ}/V^{\circ}})
  & \overset{\cup}{\longrightarrow} & H^1(N_{S^{\circ}/V^{\circ}}) \\
\bigcup && \bigcup && \bigcup \\
 H^1(\mathcal O_S(2E)) & \times & H^0(N_{S/V}(E))
  & \overset{\cup}{\longrightarrow} & H^1(N_{S/V}(3E)),
\end{array}
$$
the image of the obstruction map $\ob$ 
over $H^0(N_{S/V}(E))$ is contained in $H^1(N_{S/V}(3E))$.
The following lemma is essential to the proof of 
Theorem~\ref{thm:open surface}.
Let $d_S$ denote the restriction of the map
$d_{S^{\circ}}$ to $H^0(S,N_{S/V}(E))$.
\begin{lem}[{\cite[Proposition~2.4~(2)]{Mukai-Nasu}}]\label{lem:extension}
 Let $\partial: H^0(N_{S/V}(E)\big{\vert}_E)
 \rightarrow H^1(\mathcal O_E(2E)) \simeq H^1(N_{E/S}(E))$ 
 be the coboundary map of the exact sequence 
 \eqref{ses:normal bundle of line}$\otimes \mathcal O_S(E)$.
 Then the diagram
 \begin{equation*}
  \begin{array}{ccc}
   H^0(S,N_{S/V}(E)) & \mapright{d_{S}}
    & H^1(S,\mathcal O_S(2E)) \\ 
   \mapdown{|_E} && \mapdown{|_E} \\
   H^0(E,N_{S/V}(E)\big{\vert}_E) & \mapright{\partial} 
    & H^1(E,\mathcal O_E(2E))
  \end{array}
 \end{equation*}
 is commutative.
\end{lem}

\paragraph{Proof of Theorem~\ref{thm:open surface}.}
It suffices to show that the image
$\ob(v)\big{\vert}_E \in H^1(N_{S/V}(3E)\big{\vert}_E)$ 
of $\ob(v) \in H^1(N_{S/V}(3E))$ is nonzero.
By the definition of $v$, we have
$v\big{\vert}_E \ne 0$ in $H^0(N_{S/V}(E)\big{\vert}_E)$.
Then the line bundle $N_{S/V}(E)\big{\vert}_E \simeq \det N_{E/V}$ 
on $E$ is trivial.
Since \eqref{ses:normal bundle of line} does not split by assumption,
we have $\partial (v\big{\vert}_E)\ne 0$.
Hence by Lemma~\ref{lem:extension}, we conclude that
$$
\ob(v)\big{\vert}_E 
=d_{S^{\circ}}(v)\big{\vert}_E \cup v \big{\vert}_E
=\partial (v\big{\vert}_E) \cup v \big{\vert}_E
\ne 0. \qed
$$

\medskip

If $E$ is a $(-1)$-$\mathbb P^1$ on $S$ with
$\det N_{E/V}\simeq \mathcal O_{\mathbb P^1}$,
then the exact sequence \eqref{ses:normal bundle of line}
does not split if and only if $N_{E/V}$ is trivial.

\begin{ex}\label{ex:open cubic}
 Let $V_n$ be a smooth del Pezzo $3$-fold of degree $n \ne 8$ and
 let $E$ be a good line on $V_n$,
 i.e., $N_{E/V_n}$ is trivial (cf.~\S\ref{threefold}).
 If $S_n$ is a smooth hyperplane section of $V_n$ containing $E$,
 then there exists an obstructed infinitesimal deformation of 
 $S^{\circ}_n := S_n\setminus E$ in $V^{\circ}_n := V_n \setminus E$ .
 Indeed, let $\varepsilon: S_n \rightarrow S_{n+1}$ 
 be the blow-down of $E$ from $S_n$.
 Since $N_{S_n/V_n}\simeq -K_{S_n}$,
 $N_{S_n/V_n}(E) \simeq \varepsilon^*(-K_{S_{n+1}})$,
 and $h^0(-K_{S_{n+1}}) > h^0(-K_{S_n})$,
 there exists a global section $v$ of $N_{S_n/V_n}(E)$,
 but not of that of $N_{S_n/V_n}$.
 Then by Theorem~\ref{thm:open surface}, 
 the first order deformation of $S^{\circ}_n$ in $V^{\circ}_n$
 determined by $v$ is obstructed.
\end{ex}

In the rest of this section,
we discuss a generalization of Theorem~\ref{thm:open surface},
which will be used in the proof of Theorem~\ref{thm:main}.
Let $E$ be a disjoint union of smooth irreducible curves 
$E_i$ $(i=1,\ldots,m)$ on $S$
such that $E_i^2< 0$ and $\det N_{E_i/V}$ is trivial.
By the same symbol $E$ we also denote the divisor $\sum_{i=1}^m E_i$ 
on $S$. We define $V^{\circ}$ and $S^{\circ}$ as above and
compute the obstruction map 
$\ob: H^0(N_{S^{\circ}/V^{\circ}}) 
\rightarrow H^1(N_{S^{\circ}/V^{\circ}})$.
Then Lemma~\ref{lem:injective} allows us to 
regard $H^1(\mathcal O_S(2E))$ and $H^1(N_{S/V}(3E))$
as subgroups of $H^1(\mathcal O_{S^{\circ}})$
and $H^1(N_{S^{\circ}/V^{\circ}})$, respectively.
Then an argument similar to
\cite[Proposition~2.4~(1)]{Mukai-Nasu}
shows that the image of $H^0(N_{S/V}(E))$ by $d_{S^{\circ}}$
is contained in $H^1(\mathcal O_S(2E))$ 
and hence its image by $\ob$ is contained in $H^1(N_{S/V}(3E))$.
Moreover, we have
$$
\ob(v + v')\big{\vert}_E=\ob(v)\big{\vert}_E
$$
in $H^1(N_{S/V}(3E)\big{\vert}_E)$
for any $v \in H^0(N_{S/V}(E))$ and any $v'\in H^0(N_{S/V})$.
Indeed it follows from the definition 
of $d_{S^{\circ}}$ (cf.~\eqref{map:d-map})
that $d_{S^{\circ}}(v')$ is contained in 
$H^1(\mathcal O_S)$ and hence
\begin{align*}
\ob(v + v') &=(d_{S^{\circ}}(v) + d_{S^{\circ}}(v')) \cup (v+v') \\
&=\ob(v)
+ \underbrace{d_{S^{\circ}}(v) \cup v'+d_{S^{\circ}}(v') 
\cup v + d_{S^{\circ}}(v') \cup v'}_{\mbox{contained in $H^1(N_{S/V}(2E))$}}.
\end{align*}
Therefore the obstruction map $\ob$ induces a map
\begin{equation}\label{map:obstruction with pole}
\overline \ob: 
H^0(N_{S/V}(E))\big/ H^0(N_{S/V})
\longrightarrow H^1(N_{S/V}(3E)\big{\vert}_E). 
\end{equation}

\begin{prop}\label{prop:non-degeneracy}
 If $H^1(N_{S/V})=0$ and the exact sequence 
 \begin{equation}\label{ses2:normal bundle of line}
  0 \longrightarrow N_{E_i/S} \longrightarrow N_{E_i/V}
   \longrightarrow N_{S/V}\big{\vert}_{E_i} \longrightarrow 0 
 \end{equation}
 does not split for every $i$,
 then $\overline \ob$ is injective.
\end{prop}

This is an immediate consequence of the next lemma.

\begin{lem}
 Under the assumption of Proposition~\ref{prop:non-degeneracy},
 $\overline \ob$ is equivalent to
 the quadratic map 
 $$
 k^m \longrightarrow k^n, \qquad (a_1,\ldots,a_m) 
 \longmapsto (a_1^2,\ldots,a_m^2,0,\ldots,0)
 $$
 of diagonal type, where $n = \dim H^1(N_{S/V}(3E)\big{\vert}_E)$.
\end{lem}
\Proof
Since $H^1(N_{S/V})=0$, the source of the map $\overline \ob$
is isomorphic to $H^0(N_{S/V}(E)\big{\vert}_E)$.
Moreover there exist global sections $v_i$ of $N_{S/V}(E_i)$ 
such that $v_i\big{\vert}_E\ne 0$ in 
$H^0(N_{S/V}(E_i)\big{\vert}_{E_i})$ for all $i$.
Since $E_i$'s are mutually disjoint, we have
$N_{S/V}(E)\big{\vert}_E 
\simeq \bigoplus_{i=1}^m N_{S/V}(E_i)\big{\vert}_{E_i}
\simeq \bigoplus_{i=1}^m \mathcal O_{E_i}$.
Then there exists a commutative diagram
$$
\begin{array}{ccccccccc}
 0 & \rightarrow & H^0(N_{S/V}) & 
  \rightarrow & H^0(N_{S/V}(E)) & \rightarrow
  & H^0(N_{S/V}(E)\big{\vert}_E) & \rightarrow & 0\\
 && \mapup{a_1} && \mapup{a_2} && \mapup{a_3} && \\
 0 & \rightarrow & \bigoplus_i H^0(N_{S/V}) & 
  \rightarrow & \bigoplus_i H^0(N_{S/V}(E_i)) & \rightarrow 
  & \bigoplus_i H^0(N_{S/V}(E_i)\big{\vert}_{E_i}) & 
  \rightarrow & 0, 
\end{array}
$$
where the two horizontal sequences are exact and
$a_i$ $(1\le i\le 3)$ are defined by addition.
Since $a_1$ and $a_3$ are surjective, so is $a_2$.
Hence every element $v \in H^0(N_{S/V}(E))$ is
written as a $k$-linear combination $\sum_{i=1}^m c_i v_i$
of $v_i \in H^0(N_{S/V}(E_i))$
and the expression is unique modulo $H^0(N_{S/V})$.
By the commutative diagram
$$
\begin{array}{ccccc}
H^1(\mathcal O_S(2E)) & \times & H^0(N_{S/V}(E)) 
& \overset{\cup}{\longrightarrow} & H^1(N_{S/V}(3E)) \\
\mapdown{\vert_E} & & \mapdown{\vert_E} && \mapdown{\vert_E} \\
\bigoplus_i H^1(\mathcal O_{E_i}(2E_i)) 
& \times & \bigoplus_i H^0(N_{S/V}(E_i)\big{\vert}_{E_i}) 
& \overset{\cup}{\longrightarrow} & 
\bigoplus_i H^1(N_{S/V}(3E_i)\big{\vert}_{E_i}),
\end{array}
$$
we have
$$\ob(v) \big{\vert}_E
=(d_{S^{\circ}}(v) \cup v)\big{\vert}_E
=d_{S^{\circ}}(v)\big{\vert}_E \cup v\big{\vert}_E
=\sum_i c_i^2 d_{S^{\circ}}(v_i) \big{\vert}_{E_i} \cup v_i\big{\vert}_{E_i}.
$$
By Lemma~\ref{lem:extension}, 
$d_{S^{\circ}}(v_i) \big{\vert}_{E_i}$ is equal to
$\partial_i (v\big{\vert}_{E_i})$ in $H^1(\mathcal O_{E_i}(2E_i))$, 
where $\partial_i$ is
the coboundary map of \eqref{ses2:normal bundle of line}.
Since \eqref{ses2:normal bundle of line} does not split by assumption,
we have $\partial_i (v\big{\vert}_{E_i})\ne 0$ and hence
$d_{S^{\circ}}(v_i) \big{\vert}_{E_i}\ne 0$ for any $i$.
As a result, 
$d_{S^{\circ}}(v_i) \big{\vert}_{E_i} \cup v_i \big{\vert}_{E_i}$ 
$(1 \le i \le m)$ form a sub-basis of $H^1(N_{S/V}(3E)\big{\vert}_E)$.
\qed

\medskip

\begin{cor}\label{cor:non-degeneracy}
Let $E_i$ ($i=1,\ldots,m$) be mutually disjoint 
$(-1)$-$\mathbb P^1$'s on $S$ such that 
$N_{E_i/V} \simeq {\mathcal O_{\mathbb P^1}}^{\oplus 2}$.
If $H^1(N_{S/V})=0$, then $\overline \ob$ is injective.
\end{cor}

\section{Obstructions to deforming curves}\label{obstruction}
Let $V$ be a smooth projective $3$-fold.
In this section we study the deformation of 
smooth curves $C$ on $V$ under the presence of smooth surface 
$S$ such that $C \subset S \subset V$.
In what follows, we use the following convention.
\begin{dfn}
\begin{enumerate}
 \item $C$ is said to be {\em stably degenerate} if
       every (small) deformation of $C$ in $V$ is contained in 
       a divisor $S' \algeeq S$ of $V$
 \item $C$ is said to be {\em $S$-normal} if 
       the restriction map \eqref{map:restriction}
       is surjective.
\end{enumerate}
\end{dfn}

\subsection{$S$-normal curves and $S$-maximal families}\label{maximal}
Let $U_S$ be an irreducible component of $\Hilb V$ 
passing through $[S]$ and let 
$$
V \times  U_S \supset \mathcal S  \stackrel{p_2}\longrightarrow  U_S
$$
be the universal family of $U_S$.
Let us denote the Hilbert scheme of smooth connected curves
in $\mathcal S$ by $\Hilb^{sc} \mathcal S$, 
which is the relative Hilbert scheme of $\mathcal S/U_S$.
$\Hilb^{sc} \mathcal S$ is regarded as an open subscheme
of the Hilbert-flag scheme of $V$ (see \cite{Kleppe85} for the
definition), which parametrizes all flat families of pairs 
$(C,S)$ of a curve $C$ and a surface $S$ in $V$ such that $C \subset S$.
The projection $p_1 : \mathcal S \rightarrow V$ induces
a natural morphism 
\begin{equation}\label{map:forgetful}
pr_1: \Hilb^{sc} \mathcal S \longrightarrow \Hilb^{sc} V, 
\end{equation}
which is the forgetful morphism $(C,S) \mapsto C$.
If $pr_1$ is surjective in a neighborhood
of $[C] \in \Hilb^{sc} V$, then $C$ is stably degenerate.

Let us denote the tangential map of $pr_1$ at $(C,S)$ by
\begin{equation}\label{map:tangential}
 \kappa_{C,S}: H^0(C,N_{C/\mathcal S}) \longrightarrow H^0(C,N_{C/V}).
\end{equation}
Then we have
\begin{lem}\label{lem:universal hilb}
 Assume that $\Hilb V$ is nonsingular at $[S]$.
 If $H^1(C,N_{C/S})=0$, then
 \begin{enumerate}
  \item $\Hilb^{sc} \mathcal S$ is nonsingular at $(C,S)$, and
  \item Each of the kernel and the cokernel of $\kappa_{C,S}$ 
	is isomorphic to that of
	the restriction map \eqref{map:restriction}.
 \end{enumerate}
\end{lem}
For the proof, we refer to \cite[Lemma~1.10]{Kleppe88} for (1)
and \cite[Lemma~3.3]{Mukai-Nasu} for (2).
We can prove (2) by using the ``fundamental exact sequence 
relating $A^i(C\subset S)$ and $H^{i-1}(N_{C/V})$''
in \cite{Kleppe88} also.

\begin{prop}\label{prop:principle}
 Assume that $\Hilb V$ is nonsingular at $[S]$ and
 $H^1(C,N_{C/S})=0$.  Then:
 \begin{enumerate}
  \item If $C$ is $S$-normal, then $C$ is stably degenerate.
  \item If \eqref{map:restriction} is an
	isomorphism, then
	$\Hilb^{sc} V$ is nonsingular at $[C]$.
 \end{enumerate}
\end{prop}
\Proof
(1) Since $C$ is $S$-normal, 
$\kappa_{C,S}$ is surjective by Lemma~\ref{lem:universal hilb}~(2).
This implies that $pr_1$ is surjective in a neighborhood of $[C]$
and hence $C$ is stably degenerate.

(2) By Lemma~\ref{lem:universal hilb}~(2),
$\Hilb^{sc} \mathcal S$ is isomorphic to $\Hilb^{sc} V$
in a neighborhood of $(C,S)$.
Because $\Hilb^{sc} \mathcal S$ is nonsingular at $(C,S)$
by Lemma~\ref{lem:universal hilb}~(1), so is $\Hilb^{sc} V$ at $[C]$.
\qed

\medskip

We recall the $S$-maximal family introduced in \cite[\S3.2]{Mukai-Nasu}.
Suppose that $\Hilb V$ is nonsingular at $[S]$ and $H^1(C,N_{C/S})=0$. 
By Lemma~\ref{lem:universal hilb}~(1),
there exists a unique irreducible component
$\mathcal W_{S,C}$ of $\Hilb^{sc} \mathcal S$ containing $(C,S)$.

\begin{dfn}
We define the {\em $S$-maximal family of curves} containing $C$ 
to be the image of $\mathcal W_{S,C}$ in $\Hilb^{sc} V$
and denote it by $W_{S,C}$.
\end{dfn}

\subsection{Deformation of curves on a del Pezzo  $3$-fold}
\label{delpezzo}
Let $V$ be a smooth del Pezzo $3$-fold with the polarization $H$,
$S$ a smooth member of $|H|$, and
$C$ a smooth connected curve on $S$.
Let $n$ denote the degree of $V$ and let $d$ and $g$ denote
the degree $(:=(C\cdot H)_V)$ and the genus of $C$, respectively.

Since $-K_V\sim 2S$, by adjunction we have
$N_{S/V}= \mathcal O_S(S) \simeq - K_S$ and 
$N_{C/S}\simeq -K_S\big{\vert}_C+K_C$. 
Since $-K_S$ is ample, we have $H^1(N_{S/V})=H^1(N_{C/S})=0$.
Hence $\Hilb V$ and $\Hilb S$ are nonsingular
at $[S]$ and $[C]$, respectively.
Hence by Proposition~\ref{prop:principle}~(1), we have
\begin{lem}\label{lem:naturally degenerate}
 If $C$ is $S$-normal, then $C$ is stably degenerate.
\end{lem}
There exists a natural exact sequence
\begin{equation}\label{ses:normal bundle of curve}
0 \longrightarrow N_{C/S}
\longrightarrow N_{C/V}
\overset{\pi_S}{\longrightarrow} 
N_{S/V}\big{\vert}_C \longrightarrow 0.
\end{equation}
Since $H^1(N_{C/S})=0$, we have 
$H^1(N_{C/V})\simeq H^1(N_{S/V}\big{\vert}_C)$.
Thus every obstruction to deforming $C$ is contained in 
the cohomology group $H^1(N_{S/V}\big{\vert}_C)$.
Since  $\chi(N_{C/V})=(-K_V\cdot C)_V=2d$, we have
\begin{lem}\label{lem:naturally smooth}
 If $H^1(N_{S/V}\big{\vert}_C)=0$, then
 $\Hilb^{sc} V$ is nonsingular of expected dimension $2d$ at $[C]$.
\end{lem}
In particular, if $C$ is rational ($g=0$) or elliptic ($g=1$),
then the $\Hilb^{sc} V$ is nonsingular at $[C]$
because $H^1(N_{S/V}\big{\vert}_C)\simeq H^1(-K_S\big{\vert}_C)=0$.

Let $W_{S,C}$ be the $S$-maximal family $W_{S,C}$ of curves containing $C$.
We compute the dimension of $W_{S,C}$.
Let $pr_1: \Hilb^{sc} \mathcal S \rightarrow \Hilb^{sc} V$
be the map \eqref{map:forgetful}.
\begin{lem}\label{lem:closed embedding}
 \begin{enumerate}
  \item	$\Hilb^{sc} \mathcal S$ is nonsingular of dimension $d+g+n$
	at $(C,S)$.
  \item If $g \ge 2$ or $d \ge n+1$, then
	$pr_1$ is a closed embedding in a neighborhood of $(C,S)$
	and $\dim W_{S,C}= d+g+n$.
 \end{enumerate}
\end{lem}
\Proof
(1) Let $\mathcal W_{S,C}$ be the irreducible component of 
$\Hilb^{sc} \mathcal S$ containing $(C,S)$. 
By the Riemann-Roch theorem on $S$, 
we have $\dim |\mathcal O_{S}(C)|=d+g-1$. 
Then $\mathcal W_{S,C}$ is birationally equivalent 
to $\mathbb P^{d+g-1}$-bundle over
an open subset of $|H| \simeq \mathbb P^{n+1}$.
Hence we have $\dim \mathcal W_{S,C}= d+g+n$.

(2) By assumption, we have
$(-K_S-C)\cdot C=2-2g<0$ or $(-K_S-C)\cdot (-K_S)=n-d<0$.
Since both $C$ and $-K_S$ are nef, 
we have $H^0(N_{S/V}(-C))\simeq H^0(-K_S-C)=0$.
By Lemma~\ref{lem:universal hilb}~(2), 
$pr_1$ is a closed embedding near $(C,S)$.
Since we have $H^0(N_{S'/V}(-C'))=0$
for every generic member $(C',S')$ of $\mathcal W_{S,C}$,
the restriction of $pr_1$ to $\mathcal W_{S,C}$
is generically an embedding.
Hence $\dim W_{S,C}=\dim \mathcal W_{S,C}$.
\qed

\medskip

 We denote by $\Hilb^{sc}_{d,g} V$ 
 the open and closed subscheme of $\Hilb^{sc} V$ 
 of curves of degree $d$ and genus $g$.
 It is known that the dimension of
 every irreducible component of $\Hilb^{sc}_{d,g} V$ is greater
 than or equal to the expected dimension $\chi(N_{C/V})=2d$.

\begin{prop}\label{prop:numerically non-degenerate}
 If $\chi(V,\mathcal I_C(S))<1$, 
 then $C$ is not stably degenerate, i.e.,
 there exists a global deformation $C'$ of $C$ in $V$
 which is not contained in any divisor $S' \algeeq S$ of $V$.
\end{prop}
\Proof
There exists an exact sequence
$0 \rightarrow \mathcal I_C(S) \rightarrow 
\mathcal O_V(S) \rightarrow \mathcal O_C(S) 
\rightarrow 0$ on $V$.
We see that $\chi(\mathcal O_C(S))=d+1-g$ and $\chi(\mathcal O_V(S))=n+2$.
Hence $\chi(V,\mathcal I_C(S))<1$ is equivalent to $g < d-n$.
Then we have $\dim W_{S,C}\le \dim \mathcal W_{S,C}=d+g+n < 2d$.
Hence there exists an irreducible component $W' \supset W_{S,C}$ 
of $\Hilb^{sc} V$ such that $\dim W'> \dim W_{S,C}$.
A general member $C'$ of $W' \setminus W_{S,C}$ is
such a deformation of $C$ in $V$.
\qed

\begin{prop}\label{prop:naturally smooth}
 If $C$ is $S$-normal, then
 $\Hilb^{sc} V$ is nonsingular at $[C]$.
\end{prop}
\Proof
We may assume that $H^1(N_{S/V}\big{\vert}_C)\simeq 
H^1(-K_S\big{\vert}_C)\ne 0$ by Lemma~\ref{lem:naturally smooth}. 
Then we have $g\ge 2$.
Hence we have $H^0(N_{S/V}(-C))=0$ by Lemma~\ref{lem:closed embedding}~(2).
Then it follows from the exact sequence
\begin{equation}\label{ses:S-normal}
0 \longrightarrow N_{S/V}(-C) \longrightarrow N_{S/V}
\longrightarrow N_{S/V}\big{\vert}_C \longrightarrow 0 
\end{equation}
that the restriction map \eqref{map:restriction} is an isomorphism.
Hence $\Hilb^{sc} V$ is nonsingular at $[C]$
by Proposition~\ref{prop:principle}~(2).
\qed

\medskip

Since $H^1(N_{S/V})=0$, \eqref{ses:S-normal} shows that
$C$ is $S$-normal if and only if
$H^1(N_{S/V}(-C))=0$. There exists an exact sequence
\begin{equation}\label{ses:quotient of ideals}
[0 \longrightarrow \mathcal I_S
\longrightarrow \mathcal I_C
\longrightarrow \mathcal O_S(-C)
\longrightarrow 0] \otimes \mathcal O_V(S)
\end{equation}
on $V$. Since $\mathcal I_S(S)\simeq \mathcal O_V$ and $V$ is del Pezzo,
we have $H^i(V,\mathcal I_S(S))=0$ for $i=1,2$. 
Hence we have an isomorphism 
\begin{equation}\label{isom:abnormality}
H^1(V,\mathcal I_C(S)) \simeq H^1(S,N_{S/V}(-C)) 
\end{equation}
by $\mathcal O_S(S)\simeq N_{S/V}$.

\subsection{Stably degenerate curves}\label{proof}

We devote this subsection to the proof of Theorem~\ref{thm:main}.
Notation is same as in the previous subsection.
The following are equivalent:
(i) $\chi(V,\mathcal I_C(S))\ge 1$,
(ii) $\chi(S,N_{S/V}(-C))\ge 0$ and
(iii) $g \ge d-n$.
Indeed we have already seen in the proof of
Proposition~\ref{prop:numerically non-degenerate}
that (i) and (iii) are equivalent.
Also (i) and (ii) are equivalent because we have
$\chi(N_{S/V}(-C))=\chi(V,\mathcal I_C(S))-1$ 
by \eqref{ses:quotient of ideals}.
Throughout this subsection, we assume one of them (and hence all).

\begin{lem}\label{lem:nonspecial}
 If $H^1(N_{S/V}\big{\vert}_C)=0$ then $C$ is $S$-normal.
\end{lem}
\Proof
It suffices to show that $H^1(N_{S/V}(-C))=0$.
Since $H^2(N_{S/V})\simeq H^2(-K_S)=0$ and
$H^1(N_{S/V}\big{\vert}_C)=0$, 
we obtain $H^2(N_{S/V}(-C))=0$ by \eqref{ses:S-normal}. 
Then by assumption, we have
$0 \le \chi(N_{S/V}(-C)) = h^0(N_{S/V}(-C)) - h^1(N_{S/V}(-C))$.
Therefore if $H^0(N_{S/V}(-C))=0$, we have then $H^1(N_{S/V}(-C))=0$.
Suppose that $H^0(N_{S/V}(-C))\ne 0$.
There exists an effective divisor $D$ on $S$
such that $N_{S/V}(-C) \simeq \mathcal O_S(D)$.
If $D=0$, then $H^1(N_{S/V}(-C))=0$. Suppose that $D\ne 0$. 
Let $h$ be a general member of $|-K_S|$.
Then $h$ is a smooth elliptic curve on $S$.
Since $-K_S$ is ample, we have
$\deg \mathcal O_S(D)\big{\vert}_{h}=D \cdot (-K_S)>0$ and hence
$H^1(\mathcal O_S(D)\big{\vert}_{h})=0$.
Since $C$ is connected, 
we obtain $H^1(D-h)\simeq H^1(-C)=0$ from
the exact sequence
$0 \rightarrow \mathcal O_S(-C) \rightarrow \mathcal O_S
\rightarrow \mathcal O_C \rightarrow 0$.
Therefore it follows from the exact sequence
$$
0 \longrightarrow \mathcal O_S(D-h)
\longrightarrow \mathcal O_S(D)
\longrightarrow \mathcal O_S(D)\big{\vert}_{h}
\longrightarrow 0
$$
that $H^1(N_{S/V}(-C))\simeq H^1(D)=0$.
\qed
\medskip

Let $E_1,\ldots,E_m$ be lines on $S$ disjoint to $C$.
We define an effective divisor $E$ on $S$ by
$E:=\sum_{i=1}^m E_i$.
If $C$ is not $S$-normal, then $E$ is responsible for the abnormality.

\begin{prop}\label{prop:lifting}
 Suppose that $C$ is not rational nor elliptic.
 \begin{enumerate}
  \item  The restriction map
	 $H^0(N_{S/V}(E)) \overset{|_C}{\longrightarrow}
	 H^0(N_{S/V}\big{\vert}_C)$
	 is an isomorphism.
  \item  $C$ is $S$-normal if and only if 
	 there exists no line $\ell$ such that 
	 $C \cap \ell = \emptyset$ (i.e.~$E=0$).
 \end{enumerate}
 
\end{prop}
\Proof 
(1) Since $N_{S/V}\simeq -K_S$, we have the assertion by
Proposition~\ref{prop:anti-canonical}.

(2) The `if' part follows from (1). We prove the `only if' part.
Suppose that there exist such lines on $S$.
Let $\varepsilon: S \rightarrow F$ be the blow-down of $E$ from $S$.
Then $F$ is also a del Pezzo surface and $\varepsilon^*(-K_F)=-K_S+E$.
Since $\deg F > \deg S$, we have $h^0(-K_F)> h^0(-K_S)$.
Hence it follows from $N_{S/V}\simeq -K_S$ that
$N_{S/V}(E)$ has more global sections than $N_{S/V}$.
Hence we have $h^0(N_{S/V}\big{\vert}_C)=h^0(N_{S/V}(E))>H^0(N_{S/V})$
by (1). Therefore $C$ is not $S$-normal.
\qed

\medskip 

Let $\kappa_{C,S}: H^0(N_{C/\mathcal S}) \rightarrow H^0(N_{C/V})$ 
denote the tangential map \eqref{map:tangential}.

\begin{prop}\label{prop:core}
Suppose that $C$ is not $S$-normal. 
If every $E_i$ is a good line on $V$,
 then the obstruction $\ob(\alpha)$ is nonzero
 for any $\alpha \in H^0(N_{C/V}) \setminus \im \kappa_{C,S}$.
\end{prop}

\medskip

\Proof
We compute the exterior component $\ob_S(\alpha)$ of $\ob(\alpha)$ 
(cf.~Definition~\ref{dfn:exterior}) instead of $\ob(\alpha)$ itself.
Since $C$ is not $S$-normal, by Lemma~\ref{lem:nonspecial},
we have $H^1(N_{S/V}\big{\vert}_C)\ne 0$.
In particular, $C$ is not rational nor elliptic.
By Proposition~\ref{prop:lifting}~(1),
there exists a global section $v$ of $N_{S/V}(E)$
whose restriction $v \big{\vert}_C \in H^0(N_{S/V}\big{\vert}_C)$ 
to $C$ coincides with $\pi_S(\alpha)$.
Then $v$ is not a global section of $N_{S/V}$.
In other words, the exterior component of $\alpha$
lifts to an ``infinitesimal deformation with a pole''
(cf.~\S\ref{with pole}).
Indeed since $\alpha$ is not contained in $\im \kappa_{C,S}$,
$\pi_S(\alpha)$ is not contained in the image of \eqref{map:restriction}
by Lemma~\ref{lem:universal hilb}~(2).

Let $S^{\circ}$ and $V^{\circ}$ respectively denote
the two complements $S\setminus E$ and $V\setminus E$ of $E$.
Then $v$ determines a global section $v^{\circ}$
of the normal sheaf $N_{S^{\circ}/V^{\circ}}$
and hence the obstruction $\ob(v^{\circ}) 
\in H^1(N_{S^{\circ}/V^{\circ}})$ to deforming $S^{\circ}$
in $V^{\circ}$.
As we saw in \S\ref{with pole}, or more precisely,
by Lemma~\ref{lem:injective}, there exists a natural injection 
$H^1(S,N_{S/V}(3E)) \hookrightarrow 
H^1(S^{\circ},N_{S^{\circ}/V^{\circ}})$
and so we identify $H^1(N_{S/V}(3E))$
with its image in $H^1(N_{S^{\circ}/V^{\circ}})$.
Then $\ob(v^{\circ})$ is contained in 
the subgroup $H^1(N_{S/V}(3E))$ of $H^1(N_{S^{\circ}/V^{\circ}})$.
Now we show that $\ob(v^{\circ})$ is nonzero.
Since $v$ is not a global section of $N_{S/V}$,
the restriction $v\big{\vert}_E$ to $E$
is nonzero global section of $N_{S/V}(E)\big{\vert}_E$.
By assumption, every component $E_i$ of $E$ is a good line on $V$
and hence $N_{E_i/V}$ is a trivial bundle on $E_i \simeq \mathbb P^1$.
Therefore by virtue of Corollary~\ref{cor:non-degeneracy},
the image of $v\big{\vert}_E$ by the reduced obstruction map
$\overline \ob: H^0(N_{S/V}(E)\big{\vert}_E)
\rightarrow H^1(N_{S/V}(3E)\big{\vert}_E)$ 
(see \eqref{map:obstruction with pole} for its definition)
is nonzero, and it is equal to the restriction 
$\ob(v^{\circ})\big{\vert}_E$ of $\ob(v^{\circ})$ to $E$.
Hence we have $\ob(v^{\circ})\ne 0$ in $H^1(S,N_{S/V}(3E))$.

Finally we show that $\ob_S(\alpha)\ne 0$ in $H^1(N_{S/V}\big{\vert}_C)$.
There exists an exact sequence 
$$
0 \longrightarrow N_{S/V}(3E-C)
\longrightarrow N_{S/V}(3E) 
\overset{|_C}\longrightarrow N_{S/V}\big{\vert}_C
\longrightarrow 0.
$$
Since $N_{S/V} \simeq -K_S$, the restriction map 
$H^1(N_{S/V}(3E)) \rightarrow H^1(N_{S/V}\big{\vert}_C)$
is injective by Lemma~\ref{lem:injective2}.
Therefore we have $\ob_S(\alpha)=\ob(v^{\circ})\big{\vert}_C\ne 0$
by Lemma~\ref{lem:reduction2}.
\qed

\medskip

\paragraph{Proof of Theorem~\ref{thm:main}.}
If $C$ is $S$-normal 
then $C$ is stably degenerate 
by Lemma~\ref{lem:naturally degenerate},
and $\Hilb^{sc} V$ is nonsingular at $[C]$ 
by Proposition~\ref{prop:naturally smooth}.
Suppose that $C$ is not $S$-normal and
let $\tilde C$ be any first order deformation of $C$ in $V$.
If $\tilde C$ is contained in the image of the map $\kappa_{(C,S)}$
(cf.~\eqref{map:tangential}),
then there exists a first order deformation $\tilde S$ of $S$
such that $\tilde S \supset \tilde C$.
Since $\Hilb^{sc} \mathcal S$ is nonsingular at $(C,S)$,
the first order deformation $(\tilde C,\tilde S)$ of $(C,S)$
lifts to a global (non-infinitesimal) deformation $(C',S')$.
If $\tilde C$ is not contained in the image of $\kappa_{(C,S)}$,
then by Proposition~\ref{prop:core}, 
$\tilde C$ does not lift to a deformation over
$\Spec k[t]/(t^3)$. 
Therefore $\Hilb^{sc} V$ is singular at $[C]$ and moreover
every small deformation of $C$ in $V$ 
is contained the $S$-maximal family $W_{S,C}$ of curves
containing $C$.
Therefore $C$ is stably degenerate.
By \eqref{isom:abnormality} $C$ is $S$-normal
if and only if $H^1(V,\mathcal I_C(S))=0$.
Hence the proof of Theorem~\ref{thm:main} is completed.
\qed

\begin{rmk}
 We give two remarks on Theorem~\ref{thm:main}.
 \begin{enumerate}
  \item Suppose that $V$ is not isomorphic to a blow-up 
	$V_7$ of $\mathbb P^3$ at a point.
	If $S \in |H|$ is general, then
	by Lemma~\ref{lem:general section}, every line on $S$ 
	is a good line on $V$. Hence every curve $C$ on $S$ 
	is stably degenerate by the theorem.
	Meanwhile there exists a non-stably degenerate 
	curve $C$ on $V_7$ which is contained in a general member 
	$S$ of $|H|$ (cf.~Proposition~\ref{prop:non-stably degenerate}).
  \item There exists no line on a del Pezzo $3$-fold 
	$V_8 \simeq \mathbb P^3$.
	Hence if $V=V_8$, 
	then the assumption of the theorem concerning 
	lines $\ell$ on $S$ such that $C \cap \ell = \emptyset$ is empty.
	In fact, if $g \ge d-8$ 
	then every curve $C$ on $V_8$ is $S$-normal
	and hence stably degenerate.
	This coincides with the previous result 
	\cite[Appendix, Proposition~4.11]{Nasu},
	which proved that every curve 
	of degree $e$ and genus $p \ge 2e-8$ in $\mathbb P^3$
	lying on a smooth quadric surface
	$Q_2 \simeq \mathbb P^1 \times \mathbb P^1$ 
	is stably degenerate.
 \end{enumerate}
\end{rmk}

The following proposition is more practical 
than Proportion~\ref{prop:core}
in showing that $\Hilb^{sc} V$ is singular at $[C]$.

\begin{prop}\label{prop:singularity criterion}
 Suppose that $C$ is not rational nor elliptic. 
 If there exists a good line 
 $\ell$ on $V$ such that $\ell \subset S$ and $C \cap \ell =\emptyset$, 
 then $\Hilb^{sc} V$ is singular at $[C]$.
\end{prop}

The proofs of Proposition~\ref{prop:core} and
Proposition~\ref{prop:singularity criterion} are very similar.
Take a global section $v \in H^0(N_{S/V}(\ell))\setminus H^0(N_{S/V})$
and put $\alpha \in H^0(N_{C/V})$ as a lift of 
$v\big{\vert}_C \in H^0(N_{S/V})$ by the surjective map
$\pi_S: H^0(N_{C/V}) \twoheadrightarrow H^0(N_{S/V}\big{\vert}_C)$.
Then it is enough to show that $\ob_S(\alpha)\ne 0$ in 
$H^1(N_{S/V}\big{\vert}_C)$ by reducing it to 
$\ob(v)\big{\vert}_\ell\ne 0$ as in the proof of 
Proposition~\ref{prop:core}. We omit the details.

The following is an analogue of Conjecture~\ref{conj:Kleppe-Ellia}
due to Kleppe and Ellia.

\begin{thm}\label{thm:generic smoothness}
Let $C$ be the curve in Theorem~\ref{thm:main}. Then:
\begin{enumerate}
 \item The $S$-maximal family $W_{S,C} \subset \Hilb^{sc} V$
       of curves containing $[C]$ is an irreducible component of 
       $(\Hilb^{sc} V)_{red}$.
 \item $\Hilb^{sc} V$ is generically smooth along $W_{S,C}$ 
       if $H^1(V,\mathcal I_C(S))=0$,
       and generically non-reduced along $W_{S,C}$ otherwise.
\end{enumerate}
\end{thm}


\Proof
(1) By definition $W_{S,C}$ is an irreducible closed subset of
$\Hilb^{sc} V$. By Theorem~\ref{thm:main},
every small deformation of $C$ 
in $V$ is contained in $W_{S,C}$. This implies that
$W_{S,C}$ is a maximal irreducible closed subset of $\Hilb^{sc} V$.

(2) Let $C'$ be a general member of $W_{S,C}$.
Then $C'$ is contained in a smooth surface $S'\sim S$ in $V$.
Since $C'$ is general, so is $S'$ in $|S|$.
Suppose that $H^1(\mathcal I_C(S))=0$.
Then since $(C',S')$ is a generalization of $(C,S)$,
we have $H^1(\mathcal I_{C'}(S'))=H^1(\mathcal I_C(S))=0$ 
by the upper semicontinuity.
Hence $\Hilb^{sc} V$ is nonsingular at $[C']$ and
hence generically smooth along $W_{S,C}$.
Suppose that $H^1(\mathcal I_C(S))\ne 0$, i.e.,
$C$ is not $S$-normal.
Then Lemma~\ref{lem:nonspecial} shows that 
$H^1(N_{S/V}\big{\vert}_C)\ne 0$ and hence $g \ge 2$.
By Proposition~\ref{prop:lifting}~(2),
there exists a line $\ell$ on $S$ such that $C \cap \ell=\emptyset$.
Since $H^1(\mathcal O_S)=0$,
the Picard group of $S$ does not change under smooth deformation
and hence $\Pic S \simeq \Pic S'$.
Since $H^1(\mathcal O_S(\ell))=0$,
the line $\ell$ is deformed to a line $\ell'$ on $S'$. 
Then we have $C' \cap \ell'=\emptyset$.
Moreover since $\ell$ is a good line, so is $\ell'$.
Hence $\Hilb^{sc} V$ is singular at $[C']$ by 
Proposition~\ref{prop:singularity criterion}.
Since $C'$ is a general member of $W_{S,C}$,
$\Hilb^{sc} V$ is everywhere singular along $W_{S,C}$
and hence generically non-reduced along $W_{S,C}$.
\qed

\section{Original motivation and examples}\label{example}
\subsection{Kleppe's conjecture}\label{kleppe}
The original motivation of the present work was to show 
the following conjecture due to Kleppe.
We denote by $\Hilb^{sc}_{d,g} \mathbb P^3$ the 
open and closed subscheme of $\Hilb^{sc} \mathbb P^3$ 
consisting of curves of degree $d$ and genus $g$.
\begin{conj}[Kleppe,\,Ellia]\label{conj:Kleppe-Ellia}
Let $W$ be a maximal irreducible closed subset of 
$\Hilb^{sc}_{d,g} \mathbb P^3$ whose general member $C$ 
is contained in a smooth cubic surface. If
$$
d\ge 14, \quad g \ge 3d -18, \quad H^1(\mathbb P^3,\mathcal I_C(3))\ne 0
\quad \mbox{and} \quad H^1(\mathbb P^3,\mathcal I_C(1))= 0,
$$ 
then $W$ is a component of $(\Hilb^{sc} \mathbb P^3)_{\red}$ 
and $\Hilb^{sc} \mathbb P^3$ is generically non-reduced along $W$. 
\end{conj}
In the original conjecture \cite[Conjecture 4]{Kleppe85} of Kleppe,
the assumption of the linearly normality of $C$
(i.e.~$H^1(\mathbb P^3,\mathcal I_C(1))= 0$) was missing. 
However Ellia \cite{Ellia} pointed out that the conjecture 
does not hold for linearly non-normal curves $C$ by a counterexample, 
and suggested restricting the conjecture to linearly normal ones.
The most crucial part to prove this conjecture is the proof of the
maximality of $W$ in $(\Hilb^{sc} \mathbb  P^3)_{\red}$.
Once we prove that 
$W$ is a component of $(\Hilb^{sc} \mathbb P^3)_{\red}$,
then the non-reducedness of $\Hilb^{sc} \mathbb P^3$ 
along $W$ naturally follows.
Therefore Conjecture~\ref{conj:Kleppe-Ellia} follows from
Conjecture~\ref{conj:Kleppe-Ellia-Nasu},
where the condition 
$\chi(\mathbb P^3,\mathcal I_C(3))\ge 1$
is equivalent to $g \ge 3d-18$.
Recently it has been proved in \cite{Nasu} that 
Conjecture~\ref{conj:Kleppe-Ellia}
is true if $h^1(\mathbb P^3, \mathcal I_C(3))=1$.
Kleppe and Ellia gave a proof for the conjecture under some other
conditions, however the whole conjecture is still open.

\subsection{Hilbert scheme of canonical curves}\label{canonical}
In this subsection we prove the following:
\begin{thm}\label{thm:canonical curves}
The Hilbert scheme $\Hilb^{sc} V$ of smooth connected curves on
a smooth del Pezzo $3$-fold $V$ has a generically non-reduced component $W$.
\end{thm}
Let $n$ and $H$ be the degree and the polarization of $V$.
The theorem for the cases $n=8$ (i.e.~$V=V_8\simeq \mathbb P^3$) 
and $n=3$ (i.e.~$V$ is a smooth cubic $3$-fold $V_3$) were already obtained 
in \cite{Mumford} and \cite{Mukai-Nasu}, respectively.
For the proof, 
we consider a canonical curve $C$ on a smooth surface $S \in |H|$
which is not $S$-normal.
Here we say that a curve $C \subset V$ is {\em canonical}
if $f^*H=K_C$, where $f: C\hookrightarrow V$ is the embedding. 
Equivalently $C$ is embedded into $V$ by a linear subsystem of $|K_C|$.
Theorem~\ref{thm:generic smoothness} gives us
the non-reduced component $W$ such that $W_{\red}=W_{S,C}$.

\paragraph{Proof of Theorem~\ref{thm:canonical curves}.}
Since $V_8 \simeq \mathbb P^3$, we may assume that $n \le 7$.
Then there exists a good line $\ell$ on $V$ by Lemma~\ref{lem:good line}.
Let $S_n \in |H|$ be a smooth del Pezzo surface containing $\ell$.
We consider the complete linear system $\Lambda:=|-2K_{S_n}+2\ell|$ on
$S_n$. Let $S_{n+1}$ be the the blow-down of $\ell$ from $S_{n}$,
which is a del Pezzo surface of degree $n+1$.
Then $\Lambda$ is the pull-back of 
$|-2K_{S_{n+1}}|\simeq \mathbb P^{3n+3}$ on  $S_{n+1}$.
Since $\Lambda$ is base point free, every general member $C$ of $\Lambda$
is a smooth connected curve of degree $d=2n+2$ and genus $g=n+2$.
Therefore we have $g=d-n$ and hence $\chi(\mathcal I_C(H))=1$. 
Then $\ell$ does not intersect $C$ by 
$(-2K_{S_n}+2\ell)\cdot \ell=2-2=0$.
Moreover $\ell$ is the only such line on $S_n$.
By Theorem~\ref{thm:generic smoothness}~(1),
$W_{S_n,C}$ is an irreducible component of $(\Hilb^{sc} V)_{\red}$.
Since $C \cap \ell=\emptyset$, $C$ is not $S_n$-normal
by Proposition~\ref{prop:lifting}~(2). Therefore 
$\Hilb^{sc} V$ is generically non-reduced along $W_{S_n,C}$
by Theorem~\ref{thm:generic smoothness}~(2).
\qed

\medskip

\begin{rmk}
\begin{enumerate}
 \item By construction, $C$ is the image of a canonical curve 
 $C' \sim -2K_{S_{n+1}}$ on $S_{n+1}$
 by the projection $S_{n+1} \cdots\rightarrow S_{n}$ from a point
 $p \in S_{n+1}$ outside $C'$.
 \item The dimension of the irreducible component $W_{S_n,C}$ is equal to 
 $d+g+n=4n+4$ by Lemma~\ref{lem:closed embedding}~(2).
 \item The tangential dimension of $\Hilb^{sc} V$ at a general point $[C]$
 of $W_{S_n,C}$ is equal to $h^0(N_{C/V})=4n+5$.
 Indeed the exact sequence
 \eqref{ses:normal bundle of curve} is
 $$
 0 \longrightarrow \mathcal O_C(2K_C) \longrightarrow N_{C/V}
 \longrightarrow \mathcal O_C(K_C) \longrightarrow 0,
 $$
 since $N_{S/V}\big{\vert}_C\simeq -K_S\big{\vert}_C \simeq K_C$.
 Hence we have
 $$
 h^0(N_{C/V})=h^0(2K_C)+h^0(K_C)=(3n+3)+(n+2)=4n+5.
 $$
\end{enumerate}
\end{rmk}

The next example shows that
the curve $C$ in Theorem~\ref{thm:main}
is not necessarily stably degenerate
if there exists a bad line $\ell$ on $S$ such that 
$C \cap \ell = \emptyset$.

Let $V_7 \subset \mathbb P^8$ 
be a smooth del Pezzo $3$-fold of degree $7$
and let $S_7 \subset V_7$ be a smooth hyperplane section.
Then there exist three lines $\ell_0,\ell_1,\ell_2$ on $S_7$
forming the configuration of Figure~\ref{fig:lines on S_7}.
Consider a general member $C$ of $\Lambda:=|-2K_{S_7}+2\ell_0|$.
Then $C$ is a smooth connected curve of degree $16$ and genus $9=16-7$
and not $S_7$-normal by $C \cap \ell_0=\emptyset$. 

\begin{prop}\label{prop:non-stably degenerate}
 Let $C$ be as above.
 Then there exists a smooth deformation $C' \subset V_7$ of $C$
 not contained in any hyperplane section.
 In other words, $C$ is not stably degenerate.
\end{prop}
\Proof
Recall that $V_7$ is isomorphic to the blow-up of $\mathbb P^3$
at a point $p$. It is realized as the projection 
of the Veronese image $V_8 \subset \mathbb P^9$ of $\mathbb P^3$
from $p \in V_8$ (cf.~\S\ref{threefold}).
Then $S_7$ is the image by the projection of
a hyperplane section $Q_2 \simeq \mathbb P^1 \times \mathbb P^1$
of $V_8$ containing $p$.
Hence we have a diagram
\begin{equation}
 \begin{array}{ccccc}
  S_7 \simeq \Bl_{\text{2pts}} \mathbb P^2 & \subset 
   & V_7 \simeq \Bl_p \mathbb P^3 & \subset & \mathbb P^8 \\
  \mapdown{} \mapup{} && \mapdown{\pi_p} \mapup{\Pi_p} && \mapup{} \\
  Q_2 \simeq \mathbb P^1 \times \mathbb P^1 & \subset 
   & V_8 \simeq \mathbb P^3 & \subset & \mathbb P^9,
 \end{array}
\end{equation}
where the down arrows (resp. the up arrows) are
the blow-up morphisms at (resp. the projections from) 
$p \in Q_2 \subset V_8 \subset \mathbb P^9$.
Let $P \simeq \mathbb P^2$ denote the exceptional divisor of $\pi_p$.
Then its intersection with $S_7$ is equal to the bad line $\ell_0$.

Since $C \cap \ell_0=\emptyset$
and $C\cdot \ell_i=4$ for each $i=1,2$, $\pi_p$ maps
$C$ isomorphically onto a curve 
of bidegree $(4,4)$ on $Q_2$.
Let $Q_2'$ be a general hyperplane section of $V_8$.
Then $Q_2'\simeq \mathbb P^1 \times \mathbb P^1$ 
is mapped isomorphically onto a surface $Q_2''$ on $V_7$ by $\Pi_p$.
Here $Q_2''$ is linearly equivalent to 
$S_7 +P$ as a divisor of $V_7$ and
contains a smooth deformation $C'$ of $C$.
Then there exists no hyperplane section of $V_7$ containing $C'$.
Suppose that there exists such a hyperplane section $S_7'$.
Then the image $\pi_p(C')$ is contained in the intersection 
of two hyperplane sections $\pi_p(S_7')$ and $Q_2'$ of $V_8$.
Hence the pull-back of $\pi_p(C')$ in $\mathbb P^3$ by the Veronese 
embedding is contained in a complete intersection of two quadrics.
This is impossible since the degree of the inverse image 
is equal to $8 > 4$.
\qed

\subsection{Hilbert scheme of curves on a cubic $3$-fold}\label{cubic}
Let $V_3$ be a smooth cubic $3$-fold.
Every smooth hyperplane section $S$ of $V_3$ 
is isomorphic to a blown-up of $\mathbb P^2$ at $6$ points.
Let $\mathcal O_S(a;b_1,\ldots,b_6)$ denote
the line bundle on $S$ associated
to a divisor $a \ell - \sum_{i=1}^6 b_i e_i$ on $S$,
where $\ell$ is the pullback of a line on $\mathbb P^2$
and $e_i$ $(1 \le i \le 6)$ are the six exceptional curves on $S$.
We have an isomorphism $\Pic S \simeq \mathbb Z^7$ which sends
the class of $\mathcal O_S(a;b_1,\ldots,b_6)$
to a $7$-tuple $(a;b_1,\ldots,b_6)$ of integers.
When the linear system 
$|\mathcal O_S(a;b_1,\ldots,b_6)|$ on $S$
contains a smooth member $C$,
we denote the $S$-maximal family $W_{S,C}$ of curves containing $C$
by $W_{(a;b_1,\ldots,b_6)}$.

\begin{ex}
Suppose that $S$ is a general hyperplane section of $V_3$
and let $W$ be one of the $S$-maximal families
\begin{align*}
&W_{(\lambda +6;\lambda+1,1,1,1,1,0)} \subset
\Hilb^{sc}_{d,2d-16}V_3 \quad (d =2\lambda+13) \quad \mbox{and}\\
&W_{(\lambda+6;\lambda+2,1,1,1,1,0)} \subset
\Hilb^{sc}_{d,\frac32 d -9}V_3 \quad (d =2\lambda +12),
\end{align*}
where $\lambda \in \mathbb Z_{\ge 0}$.
Let $C$ be a general member of $W$. 
Then the genus of $C$ is greater than or equal to $d-3$
and hence $\chi(\mathcal I_C(S))\ge 1$.
Furthermore $e_6$ is the only line on $S$ such that 
$C \cap S =\emptyset$.
Since $S$ is general, $e_6$ is a good line on $V_3$ 
by Lemma~\ref{lem:general section}.
By Theorem~\ref{thm:generic smoothness},
$W$ is an irreducible component of $(\Hilb^{sc} V_3)_{\red}$
and $\Hilb^{sc} V_3$ is generically non-reduced along $W$.
Thus $\Hilb^{sc} V_3$ has infinitely many 
generically non-reduced components.
It was shown in \cite[Theorem 1.3]{Mukai-Nasu} that
for many uniruled $3$-folds $V$
the Hilbert scheme $\Hilb^{sc} V$ 
has infinitely many generically non-reduced components.
\end{ex}

\bigskip
{\sc Research Institute for Mathematical Sciences, Kyoto University,
Kyoto 606-8502, Japan}

{\it E-mail address}: nasu@kurims.kyoto-u.ac.jp

\end{document}